\documentclass[titlepage,12pt]{article}
\usepackage{amsfonts,latexsym,pstricks,multido}

\newtheorem{thm}{Theorem}[section]
\newtheorem{lem}[thm]{Lemma}
\newtheorem{defn}[thm]{Definition}

\newtheorem{prop}[thm]{Proposition}
\newtheorem{pbm}[thm]{Problem}
\newtheorem{open}[thm]{Open Problem}
\newtheorem{ex}{Example}[section]
\newtheorem{rmk}[thm]{Remark}

\newcommand{\qed}{\hspace{1em}{$\Box$}\bigskip}
\newcommand{\prf}{{\sc Proof.}$\;$}
\newcommand{\iter}{\mathrm{Iter}}
\newcommand{\Int}[2]{\mathrm{Int}_{#2}(#1)}
\newcommand{\clos}{\mathrm{Clos}}
\newcommand{\AH}{\widetilde{H}}
\newcommand{\CH}{\check{H}}

\newcommand{\usc}{usc}
\newcommand{\bw}{Bdr-w}
\newcommand{\bs}{Bdr-s}
\newcommand{\conn}{Conn}

\newcommand{\ep}{\varepsilon}

\def\R{{\mathbb R}}
\def\N{\mathbb{N}}
\def\z{\mathbb{Z}}
\def\C{\mathbb{C}}

\title{
How many times can a function be iterated?}
\author{Massimo Gobbino\vspace{2mm}\\
{\normalsize Universit\`a degli Studi di Pisa} \\{\normalsize
Dipartimento di Matematica Applicata ``Ulisse Dini''}\\
{\normalsize Via Filippo Buonarroti 1c, 56127
  PISA, Italy}\\
{\normalsize e-mail: \texttt{m.gobbino@dma.unipi.it}}\and
Robert Samuel Simon\vspace{2mm}\\
{\normalsize London School of Economics} \\{\normalsize
Department of Mathematics}\\
{\normalsize Houghton Street, WC2A 2AE, London, United Kingdom}\\
{\normalsize e-mail: \texttt{R.S.Simon@lse.ac.uk}}
}

\date{}
\begin{document}
\maketitle
\begin{abstract}
     Let $C$ be a closed subset of a topological space $X$, and let
     $f:C\to X$.  Let us assume that $f$ is continuous and $f(x)\in C$
     for every $x\in\partial C$.

     How many times can one iterate $f$?

     This paper provides estimates on the number of iterations and
     examples of their optimality.  In particular we show how some
     topological properties of $f$, $C$, $X$ are related to the maximal
     number of iterations, both in the case of functions and in the
     more general case of set-valued maps.

     We also show how this problem is related to the existence of
     equilibria for stochastic games.

\vspace{1cm}

\noindent{\bf Mathematics Subject Classification 2000 (MSC2000):}
54H20, 37B99, 55N05.

\vspace{1cm} \noindent{\bf Key words:} Viability Theory, Dynamic
Systems, iteration of functions, Point-Set Topology, Cech-Alexander
Cohomology.
\end{abstract}

\newpage

\section{Introduction}

Let $C$ be a closed subset of a topological space $X$, and let $f:C\to
X$.  We investigate the existence of finite or infinite sequences
(orbits) $\{x_{i}\}_{i\in I}$ in $X$, where $I=\{0,1,\ldots,n\}$ or $I=\N$,
such that $x_{i}=f(x_{i-1})$ for every $i\in I$ with $i\geq 1$.

At this level of generality there is of course no reason for such a
sequence to exist with $n>1$. For this reason we assume two
conditions on $f$:
\begin{itemize}
     \item  $f$ is continuous;

     \item  $f$ maps $\partial C$ back to $C$, namely $f(x)\in C$ for
     every $x\in\partial C$.
\end{itemize}

For want of a better term, we call this topic ``Discrete Viability
Theory''.  Surprisingly this problem seems to be quite new.  Up to now
indeed we have found little related literature, although this topic
seems to come close to different areas. Let us mention some of them.
\begin{itemize}
         \item \emph{Conventional Viability Theory}.  This theory, for
         which we refer the reader to \textsc{J.\ P.\ Aubin}~\cite{aubin},
         considers \emph{continuous-time} dynamic processes with some
         control mechanism.  The main problem is finding conditions under
         which these processes stay within a given set $C$.  As in our
         problem, these conditions often involve the behavior of the flow
         at the boundary of $C$.  Unfortunately there are relatively few
         theorems for \emph{discrete-time} models: one example is Theorem
         3.7.11 of \cite{aubin} where $X=\R^{n}$, $C$ is a convex subset,
         $f$ is a multi-valued map with convex images, and the existence of
         a fixed point is proven.

         \item \emph{Fixed point theorems}.  If $f$ has a fixed point $x\in
         C$, then we can clearly iterate $f$ infinitely many times starting
         from $x$.  So the most interesting case is when $f$ has no fixed
         point.  Let us assume however that an infinite orbit exists.  Then
         under general assumptions the $\omega$-limit of this orbit is a
         closed $f$-invariant set, hence a fixed point of $f$ as a
         function acting on the space of closed sets.  In this way the
         existence of an infinite orbit is reduced to a fixed point
         problem.  Unfortunately up to now this approach didn't work
         because it's difficult to find topological obstructions in the
         space of closed subsets.

	\item \emph{Dynamical systems}.  On the one hand our problem can
	be considered as a problem in discrete-time dynamical systems or
	in topological dynamics.  On the other hand, to our knowledge in
	all the literature the iterations are always well defined for the
	trivial reason that $C=X$, and the main questions concern their
	asymptotic behavior.  Continuous-time dynamical systems are weakly
	related to our topic if we add the assumption that $f$ is
	homotopic to the identity.  In this case indeed $f(x)$ could be
	interpreted as the position at time $t=1$ of a continuous
	trajectory which starts from $x$ at time $t=0$.  If this is the
	case we could apply the classical tools for the study of flows,
	such as for example the Conley Index Theory (see \cite{MM}).
	However in general there is no flow which connects $x$ and $f(x)$,
	and for this reason we dismissed this approach as hopeless. 
	One could argue that the construction known as ``suspension of a
	map'' (or ``the mapping torus'') describes how to turn a map on a
	space $X$ into a flow on a \emph{different} topological space 
	$Y$, but this could at most provide trajectories in $Y$.

	\item \emph{Game Theory}.  One of the main problems in game
	theory is the existence of equilibria.  Classical results in
	this field are usually proved by means of fixed point
	theorems (see for example the celebrated result by
	\textsc{J.\ Nash}~\cite{nash}).  In Section~\ref{sec:games} we
	show that there are special stochastic games called quitting
	games for which the existence of approximate equilibria is
	equivalent to the existence of suitable \emph{non-stationary}
	orbits for some multi-valued functions.  Unfortunately the
	counterexamples we present in this paper, one of which
	(Example~\ref{ex:corr-main}) was inspired by a game theoretic
	context, show that this equivalence doesn't lead to a simple
	proof of existence of equilibria.
\end{itemize}

Let us come now to an explicit example of a question.  Let us
consider the case where $X=\R$, and $C$ is an interval.  Our
assumption on the behavior of $f$ at the boundary implies that
$f(0)\geq 0$ and $f(1)\leq 1$, hence $f$ admits a fixed point because
of the Intermediate Value Theorem. Therefore the simplest nontrivial
problem in Discrete Viability Theory is probably the following one.

\begin{pbm}\label{Pbm:main}
   Let $C$ be a compact and connected subset of $X= \R^{n}$. Let
     $f:C\to X= \R^{n}$ be a continuous function such that $f(x)\in C$ for
     every $x\in\partial C$.

     How many times can we iterate $f$ starting from a suitable $x\in
     C$?
\end{pbm}

What motivates us the most is the search for some theorem that implies
the existence of an infinite orbit for functions that neither have
fixed points nor $f(C)\subseteq C$.  When we started attacking this
problem we were rather optimistic about the existence of such a
theorem under general assumptions on $f$, $C$, $X$, and all the
colleagues we contacted in that period shared our optimism.  The first
two iterations are indeed given for free, and a simple connectedness
argument provides two more iterations.  At a first glance it seemed
also possible to reiterate the argument (see Remark~\ref{rmk:wrong})
assuming only the connectedness of $X$ and $C$.

Our optimism decreased when \textsc{T.\ Wiandt} \cite{TW} showed us a
simple situation (see Example \ref{ex:XC-conn}) where $X$ and $C$ are
compact and connected but only four iterations are possible.  That
example showed us that further requirements on $f$, $C$, $X$ were
needed in order to perform further iterations of $f$.  In order to
rule out the situation of Example~\ref{ex:XC-conn} we worked in two
different directions: either by asking that $f$ is homotopic to the
identity in a suitable sense (see Problem~\ref{Pbm:map-homot}), or by
requiring $X$ to be simply connected, since in that example $X$ is the
unit circle $S^{1}$.  In both cases we succeeded in proving the
existence of a fifth iteration (Theorem~\ref{thm:dvt-maps},
Theorem~\ref{thm:dvt-maps-homot} and Theorem~\ref{thm:dvt-corr}).
However, the argument is more involved, and surprisingly more or less
the same despite of the different additional assumptions.

The little optimism left become pessimism when we found
Example~\ref{ex:stargate}, where $f$, $C$, $X=\R^2$ are as in
Problem~\ref{Pbm:main}, $f$ is homotopic to the identity in the
suitable sense, and nevertheless only six iterations can be computed.

In any case we are not sure that this is the end of the story, because
probably further topological requirements on $C$ can provide more
iterations (see Section~\ref{sec:open}).
The space  $\R^2$ is very
restrictive concerning the
topological options for a subset $C$ that is the closure of an open
set. It is conceivable that some additional
conditions
that would imply  the existence of an infinite orbit  (but not necessarily a
fixed point)  in higher dimensions
would   imply the existence of a fixed point      in the  $X=\R^2$
context.

In this paper we present some lower bounds for the number of
iterations in terms of the topological properties of $f$, $C$, $X$,
and we show their optimality with some examples.  In order to give a
complete theory we work both with functions and with set-valued maps
(a good reference on iterating set-valued maps is \cite{Mc}).  Many
parts of the theory are similar in both cases, but there are also some
remarkable differences (see Remark~\ref{rmk:corr}).

This paper is organized as follows.  In Section \ref{sec:statements}
we state the questions and our results.  In Section \ref{sec:proofs}
we prove the results.  In Section \ref{sec:examples} we present some
examples showing the optimality of our estimates.  In
Section~\ref{sec:games} we present the connections with game theory
which motivated this study.  In Section \ref{sec:open} we state some
open problems.

\section{Statements}\label{sec:statements}

Throughout this paper, unless otherwise stated, $X$ denotes a
topological space.  Given $A\subseteq X$, $\Int{A}{}$ denotes the set
of interior points, $\clos(A)$ the closure, $\partial A$ the boundary
of $A$ in $X$.  We recall that $X$ is said to be locally connected if
every $x\in X$ has a fundamental system of connected neighborhoods.

Every $Y\subseteq X$ may be regarded as a topological space itself,
with the topology inherited as a subset of $X$.  If now $A\subseteq
Y$, then $\Int{A}{Y}$, $\clos_{Y}(A)$, $\partial_{Y}A$ denote,
respectively, the set of interior points, the closure, and the boundary
of $A$ relative to the topological space $Y$.

We say that $Y$ satisfies the fixed point property if every continuous
function $g:Y\to Y$ has a fixed point.  For example, any nonempty
compact convex subset of $\R^{n}$ has the fixed point property 
because of Brouwer's fixed point theorem.

In this paper we make a mild use of Cech-Alexander cohomology, in the
sense that in some statements we assume that $\CH^{1}(X)=0$, namely
that the first Cech-Alexander cohomology group (with $\z$ as
coefficient group, just to fix the ideas) is trivial.  For readers
which are not familiar with this cohomology theory, in
Lemma~\ref{lemma:cech} we show that for reasonable spaces (\emph{e.g.}
paracompact Hausdorff spaces) this assumption implies the following:
``for every open set $A\subseteq X$, if $A$ and $X\setminus A$ are
connected, then $\partial A$ is connected''.  This last property is
what we use in this paper.  We recall also that a simple case in which
$\CH^{1}(X)=0$ is when $X$ is locally contractible and simply
connected.  Good references for Cech-Alexander cohomology are Chapter
3 of \cite{hu} and Chapter 6 of \cite{spanier}.

\subsection{DVT for functions}

The following is the main question in what we called Discrete
Viability Theory.

\begin{pbm}\label{Pbm:map}

         Let $X$ be a topological space, and let $C\subseteq X$ be a
         nonempty closed subset.  Let $f:C\to X$ be a continuous function
         such that $f(x)\in C$ for every $x\in\partial C$.

     How many times can we iterate $f$ starting from a suitable $x\in
     C$?
\end{pbm}

In the following problem we strengthen the assumptions on $f$ by
asking that $f$ is homotopically equivalent to the identity map on $C$
by a homotopy whose intermediate maps also send $\partial C$ back to
the set $C$.

\begin{pbm}\label{Pbm:map-homot}
     Let $X$, $C$, and $f$ be as in Problem \ref{Pbm:map}. Let us assume
     that there exists a function $\Phi:C\times [0,1]\to X$ such that
     \begin{itemize}
         \item  $\Phi(x,0)=x$ for every $x\in C$;

         \item  $\Phi(x,1)=f(x)$ for every $x\in C$;

         \item  $\Phi(x,t)\in C$ for every $x\in\partial C$ and every
         $t\in[0,1]$.
     \end{itemize}

     How many times can we iterate $f$ starting from a suitable $x\in
     C$?
\end{pbm}

In order to better investigate these problems, we introduce some
notations.

\begin{defn}\label{defn:map}
         Let $X$, $C$, and $f$ be as in Problem \ref{Pbm:map}.  We
         recursively define a sequence $\{C_{n}\}_{n\in\N}$ of subsets of
         $X$ by
         $$C_{0}:=X,
         \hspace{4em}
         C_{n+1}:=\{x\in C:\ f(x)\in C_{n}\}.$$

         Then we set
         $$A_{n} := C_{n}\setminus C_{n+1};$$
         $$\iter(f,C,X) := \sup\{n\in\N:\
                 C_{n}\neq\emptyset\}\in\N\cup\{+\infty\}.$$
\end{defn}

The following proposition clarifies the set-theoretic properties of
the notions we have just introduced (proofs are trivial).

\begin{prop}\label{prop:dvt-st-maps}
         Let $X$ be a set, let $C\subseteq X$ be a nonempty subset, and let
         $f:C\to X$ be any function.

     Then the notions introduced in Definition \ref{defn:map} fulfil
         the following properties:

     \begin{enumerate}
         \renewcommand{\labelenumi}{(\arabic{enumi})}

                 \item $\iter(f,C,X)$ is the maximal length of a sequence
                 $x_{0},\ldots,x_{n}$ such that $x_{i}=
                 f(x_{i-1})$ for every $i=1,\ldots,n$;

         \item $C_{n+1}\subseteq C_{n}$ for every $n\in\N$;

         \item  \label{if-equal}
         if $C_{n+1}=C_{n}$ for some $n\in\N$, then
         $C_{m}=C_{n}$ for every $m\geq n$;

                 \item \label{stm:st-fcc}
                 if $x\in C_{n+1}$ then $f(x)\in C_{n}$;

                 \item \label{stm:st-faa}
                 if $x\in A_{n+1}$ then $f(x)\in A_{n}$;

                 \item
                 if $\iter(f,C,X)=k<+\infty$, then $A_{i}\neq\emptyset$ if and
                 only if $i\leq k$.
     \end{enumerate}
\end{prop}

We state now the topological properties of the sets $A_{n}$ and
$C_{n}$.

\begin{prop}\label{prop:dvt-maps}
     Let $X$, $C$, and $f$ be as in Problem \ref{Pbm:map}.

         Then for every $n\in\N$ we have that (for simplicity we use 
$\partial_{n}$
         instead of $\partial_{C_{n}}$ to denote boundaries relative to
         $C_{n}$)
     \begin{enumerate}
         \renewcommand{\labelenumi}{(\arabic{enumi})}

         \item $C_{n}$ is a closed subset of $X$;

         \item \label{bndr-bndr}
 
$f\left(\partial_{n+1}C_{n+2}\right)\subseteq\partial_{n}C_{n+1}$;

         \item \label{bndr-inside}
         $\partial_{n}C_{n+1}\subseteq C_{n+2}$;

         \item \label{ak-closed}
         $A_{n}\cup C_{n+2}$ is a closed set.

     \end{enumerate}
\end{prop}

\begin{rmk}\label{rmk:wrong}
     \begin{em}
                 As a consequence of Proposition \ref{prop:dvt-st-maps} and
                 Proposition \ref{prop:dvt-maps}, by restricting the domain and
                 the codomain, we can regard $f$ as a function $f:C_{n+1}\to
                 C_{n}$, and this restriction satisfies
                 $f(\partial_{n}C_{n+1})\subseteq C_{n+1}$.  Therefore, if
                 $f:C\to X$ satisfies the assumptions of Problem \ref{Pbm:map},
                 then $f:C_{n+1}\to C_{n}$ satisfies the same assumptions for
                 every $n\in\N$, and $\iter(f,C,X)=n+\iter(f,C_{n+1},C_{n})$.

                 If we know \emph{a priori} that $C_{n}$ is connected for every
                 $n\in\N$, this leads to an inductive proof that
                 $C_{n}\neq\emptyset$ for every $n\in\N$.  But we can find no
                 non-trivial condition that forces this to hold, and simple
                 examples can be given where infinite orbits (and also fixed
                 points) exist and $C_0$ and $C_1$ are  connected, 
     but $C_n$ in not connected for all $n\geq 2$.

                 If $C_{n}$ is not connected it may happen that $C_{n+1}$ is
                 the union of some connected components of $C_{n}$: in this
                 case $f$ can map $C_{n+1}$ into the remaining connected
                 components of $C_{n}$, causing $C_{n+2}$ to be empty (see the
                 examples in Section~\ref{sec:examples}).

                 This points out once more the importance of \emph{relative}
                 boundaries in Proposition \ref{prop:dvt-maps}: boundaries are
                 always defined relative to something, and that something can
                 change at each step.
     \end{em}
\end{rmk}

The following  result provides our estimates on the number of
iterations for Problem~\ref{Pbm:map}.

\begin{thm}\label{thm:dvt-maps}
     Let $X$, $C$, and $f$ be as in Problem \ref{Pbm:map}, and let
     $\iter(f,C,X)$ be as in Definition \ref{defn:map}.
     Then we have the following estimates.
         \begin{enumerate}
             \renewcommand{\labelenumi}{(\arabic{enumi})}

             \item \label{stm:dvt2}
             If $\partial C\neq\emptyset$ then $\iter(f,C,X)\geq 2$.

             \item \label{stm:dvt3}
             If $X$ is connected, then $\iter(f,C,X)\geq 3$.

             \item \label{stm:dvt4}
             If $X$ is connected, and $C$ is connected, then
             $\iter(f,C,X)\geq 4$.

              \item \label{stm:dvt6}
                  Let us assume that $C$ is connected, and that $X$ is a
                  paracompact Hausdorff space which is connected, locally
                  connected and satisfies $\CH^{1}(X)=0$.

                  Then $\iter(f,C,X)\geq 5$.

                 \item \label{stm:dvt-infty}
                 If $\partial C$ is a retract of $X\setminus \Int{C}{}$, and
                 $C$ satisfies the fixed point property, then there exists
                 $x\in C$ such that $f(x)=x$.  In particular
                 $\iter(f,C,X)=+\infty$.
\end{enumerate}

\end{thm}

Under the assumptions of Problem~\ref{Pbm:map-homot} we have the
following result (note that there are no topological requirements on
$C$ and $X$).

\begin{thm}\label{thm:dvt-maps-homot}
         Let $X$, $C$, and $f$ be as in Problem ~\ref{Pbm:map-homot}, and let
     $\iter(f,C,X)$ be as in Definition \ref{defn:map}.

         Then we have that $\iter(f,C,X)\geq 5$.
\end{thm}

Some examples in Section~\ref{sec:examples} show the optimality of
these estimates.

\subsection{DVT for set-valued maps}

In this section we extend some parts of the theory from functions to
set-valued maps.  Let us begin with some notations and definitions.

Let $X$ be a topological space, let $C\subseteq X$ be a closed subset,
and let $\mathcal{P}_{\star}(X)$ be the set of nonempty subsets of
$X$.  A set-valued map on $C$ with values in $X$ is any map
$f:C\to\mathcal{P}_{\star}(X)$.

The first thing we need is some continuity of $f$.  There are several
notions of continuity for set-valued maps, and all of them are
equivalent to standard continuity in the case of single-valued maps.
The notion we use in this paper is usually referred in the literature
as upper semicontinuity, and it is defined as follows.

\begin{list}{}{\leftmargin 3em \labelwidth 3em}
         \item[(\usc)] A map $f:C\to\mathcal{P}_{\star}(X)$ is upper
         semicontinuous if for every open set $U\subseteq X$ we have that
         $\{x\in C:\ f(x)\subseteq U\}$ is an open subset of $C$.
\end{list}

Then we need to control the behavior of $f$ at $\partial C$.  The
assumption in Problem \ref{Pbm:map} can be extended to set-valued maps
in a weak and in a strong sense (equivalent if $f$ is single-valued),
as follows.
\begin{list}{}{\leftmargin 5em \labelwidth 5em}
     \item[(\bw)] For every $x\in\partial C$ we have that $f(x)\cap
     C\neq\emptyset$.

     \item[(\bs)] For every $x\in\partial C$ we have that
     $f(x)\subseteq C$.

\end{list}

Finally, simple examples (see Example \ref{ex:trivial-corr}) show
that nothing but the trivial iterations can be expected without
connectedness assumptions on the images.  For this reason, we often
need the following property (trivially satisfied by functions).
\begin{list}{}{\leftmargin 5em \labelwidth 5em}
     \item[(\conn)] For every $x\in C$ we have that $f(x)$ is connected.

\end{list}

We can now state the main question in Discrete Viability Theory for
set-valued maps.

\begin{pbm}\label{Pbm:corr}
         Let $X$ be a topological space, let $C\subseteq X$ be a nonempty
         closed subset, and let $f:C\to\mathcal{P}_{\star}(X)$ be a
         set-valued map satisfying (\usc), (\bw) or (\bs), and (\conn).

     How many times can we iterate $f$ starting from a suitable $x\in
     C$?
\end{pbm}

In order to study this problem, in analogy with the case of functions
we consider the sequence of sets $\{C_{n}\}_{n\in\N}$ recursively
defined by $$C_{0}:=X, \hspace{4em}
C_{n+1}:=\{x\in C:\ f(x)\cap C_{n}\neq\emptyset\},$$
and then we define $A_{n}$ and $\iter(f,C,X)$ as in Definition
\ref{defn:map}.

The set-theoretic properties of these notions are analogous to the
case of functions. We sum them up in the following Proposition.

\begin{prop}\label{prop:dvt-st-corr}
     Let $X$ be a set, let $C\subseteq X$ be a nonempty subset, and let
     $f:C\to\mathcal{P}_{\star}(X)$.

         Then statements (2), (3), (6) of Proposition
         \ref{prop:dvt-st-maps} hold true without changes.  Moreover,
         statements (1), (\ref{stm:st-fcc}), (\ref{stm:st-faa}) of
         Proposition \ref{prop:dvt-st-maps} hold true in the following
         modified form:

     \begin{enumerate}
         \renewcommand{\labelenumi}{(\arabic{enumi})}

                 \item[(1$'$)] $\iter(f,C,X)$ is the maximal length of 
a sequence
                 $x_{0},\ldots,x_{n}$ such that  $x_{i}\in
                 f(x_{i-1})$ for every $i=1,\ldots,n$;

                 \item[(\ref{stm:st-fcc}$'$)] if $x\in C_{n+1}$ then $f(x)\cap
                 C_{n}\neq\emptyset$;

                 \item[(\ref{stm:st-faa}$'$)] if $x\in A_{n+1}$ then
                 $f(x)\subseteq A_{0}\cup\ldots\cup A_{n}$ and $f(x)\cap
                 A_{n}\neq\emptyset$.
     \end{enumerate}
\end{prop}

The topological properties of the sets $A_{n}$ and $C_{n}$ are
analogous to the case of functions only for small values of $n$, as
stated in the following Proposition.

\begin{prop}\label{prop:dvt-corr}
     Let $X$ be a topological space, let $C\subseteq X$ be a closed
     subset, and let $f:C\to\mathcal{P}_{\star}(X)$ be a set-valued map
     satisfying (\usc), (\bw) and (\conn).

         Then (we use $\partial_{n}$ instead of $\partial_{C_{n}}$ to
         denote boundaries relative to $C_{n}$)
     \begin{enumerate}
     \renewcommand{\labelenumi}{(\arabic{enumi})}

         \item \label{stm:c-cnclosed}
         $C_{n}$ is a closed subset of $X$ for every $n\in\N$;

     \item \label{stm:c-a01open}
         $A_{0}$ and $A_{1}$ are open subsets of $X$;

     \item  \label{stm:c-bdr013}
         $\partial_{1}C_{2}\subseteq C_{3}$;

         \item  \label{stm:c-a02open}
         $A_{1}\cup C_{3}$ is a closed subset of $X$, hence
         $A_{0}\cup A_{2}$ is an open subset of $X$;

     \item  \label{stm:c-bdr124}
         $\partial_{2}C_{3}\subseteq C_{4}$;

         \item  \label{stm:c-a013open}
         $A_{2}\cup C_{4}$ is a closed subset of $X$, hence
         $A_{0}\cup A_{1}\cup A_{3}$ is an open subset of $X$;

         \item  \label{stm:c-fbdr34}
         if $x\in\partial_{3}C_{4}\cap A_{4}$ then $f(x)\cap
         A_{2}\neq\emptyset$ and $f(x)\cap A_{3}\neq\emptyset$.

     \end{enumerate}
\end{prop}

The following result is the counterpart of Theorem \ref{thm:dvt-maps}
for set-valued maps.

\begin{thm}\label{thm:dvt-corr}
         Let $X$ be a topological space, let $C\subseteq X$ be a closed
         subset, let $f:C\to\mathcal{P}_{\star}(X)$ be a set-valued map,
         and let $\iter(f,C,X)$ be as in Definition \ref{defn:map}.

         Then we have the following estimates.
         \begin{enumerate}
             \renewcommand{\labelenumi}{(\arabic{enumi})}

             \item \label{stm:dvtc2}
                         If $\partial C\neq\emptyset$ and $f$ 
satisfies (\bw) then
                         $\iter(f,C,X)\geq 2$.

             \item \label{stm:dvtc3}
                         Let us assume that $X$ is connected, and $f$ satisfies
                         (\usc), (\bw), and (\conn).  Then $\iter(f,C,X)\geq 3$.

             \item \label{stm:dvtc4}
                         Let us assume that $X$ is connected, $C$ is 
connected, and
                         $f$ satisfies (\usc), (\bw), and (\conn).  Then
                         $\iter(f,C,X)\geq 4$.

             \item \label{stm:dvtc5}
                         Let us assume that
                         \begin{itemize}
                                 \item $X$ is a paracompact Hausdorff 
space which is
                                 connected, locally connected and satisfies
                                 $\CH^{1}(X)=0$;

                                 \item  $C$ is connected;

                                 \item $f$ satisfies (\usc), (\bs), and (\conn).
                         \end{itemize}
             Then $\iter(f,C,X)\geq 5$.

         \end{enumerate}

\end{thm}

The optimality of these estimates follows from the optimality of the
corresponding estimates for functions.

\begin{rmk}\label{rmk:corr}
\begin{em}
                 Example \ref{ex:corr-main} shows that statement
                 (\ref{stm:dvtc4}) is the best one can expect under assumption
                 (\bw) (note that we assumed (\bs) in statement
                 (\ref{stm:dvtc5})).  In that example indeed $X$ is $\R^{2}$,
                 $C$ is a contractible compact set which satisfies the fixed
                 point property for functions, and all images of $f$ are convex
                 sets.
\end{em}
\end{rmk}

\section{Proofs}\label{sec:proofs}

\subsection{Topological lemmata}

The five lemmata we collect in this section are the technical core of
this paper.

The first one is standard point-set topology.  The statements may seem
trivial: nevertheless, at least (\ref{lemma:b-b}),
(\ref{lemma:b-b-c}), and (\ref{lemma:bdr-union}) are false without
local connectedness assumptions.

\begin{lem}\label{lemma:top}
         Let $Y$ be a locally connected topological space.

         Then the following implications are true.
     \begin{enumerate}
                 \renewcommand{\labelenumi}{(\arabic{enumi})}
                 \item \label{lemma:b-b}
                 If $V\subseteq Y$ is any subset, and $V'$ is a connected
                 component of $V$, then $\partial V'\subseteq\partial V$.

                 \item \label{lemma:b-b-c}
                 If $V\subseteq Y$ is closed, and $V'$ is a connected
                 component of $V$, then $\partial V'=V'\cap\partial V$.

                 \item \label{lemma:bdr-union}
                 For every family $\{A_{i}\}_{i\in I}$ of subsets of $Y$ we
                 have that
                 $$\partial\left(\bigcup_{i\in I}A_{i}\right)\subseteq
                 \clos\left(\bigcup_{i\in I}\partial A_{i}\right).$$

                 \item \label{lemma:extr-conn} Let us assume that $Y$ is
                 connected, $A\subseteq Y$ is an open subset such that
                 $Y\setminus A$ is connected, and $A'$ is a connected component
                 of $A$. Then $Y\setminus A'$ is connected.

     \end{enumerate}
\end{lem}
\prf

\emph{Statement (1)}.  Let $x\in\partial V'$.  Then
$x\in\clos(V')\subseteq\clos(V)$, hence either $x\in\partial V$ or
$x\in\Int{V}{}$.  Assume by contradiction that $x\in\Int{V}{}$.  Since
$Y$ is locally connected there exists a connected neighborhood $U$ of
$x$ contained in $V$.  Since $U$ is connected it is necessarily
contained in $V'$, but this implies that $x\in\Int{V'}{}$ and
contradicts the assumption that $x\in\partial V'$.  \medskip

\emph{Statement (2)}.  We have that $\partial V'\subseteq V'$
because $V'$ is closed, and $\partial V'\subseteq\partial V$ because
of the statement (\ref{lemma:b-b}).  The opposite inclusion
$V'\cap\partial V\subseteq\partial V'$ is trivial (it holds true also
without the local connectedness of $Y$ or the closedness of $V$).
\medskip

\emph{Statement (3)}.  Let $x$ be a point in the boundary of the
union, and let $U$ be any connected neighborhood of $x$.  By
assumption there exists $i_{0}\in I$ such that $U\cap
A_{i_{0}}\neq\emptyset$ and $U\setminus A_{i_{0}}\neq\emptyset$.  By
the connectedness of $U$ this implies that $U\cap\partial
A_{i_{0}}\neq\emptyset$.  Since $x$ has a fundamental system of
connected neighborhoods, this is enough to conclude that $x$ belongs
to the closure of the union of the boundaries.  \medskip

\emph{Statement (4)}.  If $A$ is connected the conclusion is
trivial.  Otherwise, let $\{A_{i}\}_{i\in I}$ be the set of connected
components of $A\setminus A'$ so that
$$Y\setminus A'=(Y\setminus A)\cup\bigcup_{i\in I}^{}A_{i}=
\bigcup_{i\in I}^{}\left[(Y\setminus A)\cup A_{i}\right].$$

Thus it is enough to show that $(Y\setminus A)\cup A_{i}$ is connected
for every $i\in I$.  Since $A_{i}$ is a nontrivial subset of the
connected space $Y$, we have that $\partial A_{i}\neq\emptyset$, hence
by statement~(\ref{lemma:b-b})
$$\emptyset\neq \partial A_{i}\subseteq\partial A=\partial(Y\setminus
A)\subseteq Y\setminus A.$$

Since $\clos(A_{i})$ is also connected it follows that $(Y\setminus
A)\cup A_{i}=(Y\setminus A)\cup \clos(A_{i})$ is the union of two
connected sets with nonempty intersection, hence it is connected.
\qed

The second lemma relates the cohomological assumption on the space to
the connectedness of the boundary of suitable subsets.  We use this
result every time we want to prove that the boundary of an open set is
connected.

\begin{lem}\label{lemma:cech}
         Let $Y$ be a paracompact Hausdorff topological space such that
         $\CH^{1}(Y)=0$. Let $A\subseteq Y$ be a connected open set such
         that $Y\setminus A$ is also connected.

         Then $\partial A$ is connected.
\end{lem}

\prf
We recall that a topological space is connected if and only if its
0-dimen\-sional reduced Alexander cohomology group (with any coefficient
group) is trivial.

Let us consider the long exact sequence of reduced Alexander cohomology
groups for the pair $(Y,\clos(A))$ (see \cite[Theorem 2.13]{hu}):
$$\ldots\longrightarrow \AH^{0}(\clos(A))\longrightarrow
\AH^{1}(Y,\clos(A))\longrightarrow
\AH^{1}(Y)\longrightarrow\ldots$$

In this sequence we have that $\AH^{0}(\clos(A))=0$ because
$\clos(A)$ is connected, and $\AH^{1}(Y)=0$.  This implies that
$\AH^{1}(Y,\clos(A))=0$.

By the strong excision property in paracompact Hausdorff spaces (see
Exercise 6B in \cite[p.  89]{hu} or Theorem 5 in \cite[p.
318]{spanier}) we can subtract $A$ to both $Y$ and $\clos(A)$ obtaining
that
$$\AH^{1}(Y\setminus A,\partial A)=\AH^{1}(Y,\clos(A))=0.$$

Thus in the long exact sequence for the pair $(Y\setminus A,\partial A)$
$$\ldots\longrightarrow \AH^{0}(Y\setminus A)\longrightarrow
\AH^{0}(\partial A)\longrightarrow
\AH^{1}(Y\setminus A,\partial A)\longrightarrow
\ldots$$
we have that $\AH^{1}(Y\setminus A,\partial A)=0$ and $\AH^{0}(Y\setminus
A)=0$ because $Y\setminus A$ is connected. It follows that $\AH^{0}(\partial
A)=0$, which is equivalent to say that $\partial A$ is connected.
\qed

The following result is used in the sequel every time we prove the
existence of a fifth iteration.  We state and prove it under the joint
hypotheses of Lemma~\ref{lemma:top} and Lemma~\ref{lemma:cech}.  We
suspect it can be true also without the local connectedness
assumption, but in that case the proof could be much more involved.
On the contrary, the cohomological assumption is likely to be
necessary.

\begin{lem}\label{lemma:absorb}
         Let $Y$ be a paracompact Hausdorff locally connected
         topological space such that
         $\CH^{1}(Y)=0$. Let $K_{1}$, $K_{2}$, $U$ be three
         subsets such that
         \begin{enumerate}
                 \renewcommand{\labelenumi}{(\roman{enumi})}
                 \item  $K_{1}\cap K_{2}=\emptyset$;

                 \item $U$ is open and $Y\setminus U$ is connected;

                 \item  $\partial U\subseteq K_{1}\cup K_{2}$;

                 \item  $\partial U\cap K_{1}$ and $\partial U\cap K_{2}$ are
                 closed sets.
         \end{enumerate}

         Then $U$ is the disjoint union of two subsets $U_{1}$ and $U_{2}$
         such that $\partial U_{1}\subseteq\clos(K_{1})$ and
         $\partial U_{2}\subseteq\clos(K_{2})$.
\end{lem}

\prf
Let $U'$ be any connected component of $U$.  By (ii) and
statement~(\ref{lemma:extr-conn}) of Lemma~\ref{lemma:top} we have
that $Y\setminus U'$ is connected, and therefore from Lemma
\ref{lemma:cech} we deduce that $\partial U'$ is connected.  Due to
statement~(\ref{lemma:b-b}) of Lemma~\ref{lemma:top} and assumption
(iii) we have that $\partial U'\subseteq\partial U\subseteq K_{1}\cup
K_{2}$.  We can therefore write
$$\partial U'=(\partial U'\cap K_{1})\cup(\partial U'\cap K_{2}).$$

By assumptions (i) and (iv), the two terms in the right hand side are
closed and disjoint, hence one of them must be empty.  This proves
that every connected component $U'$ of $U$ satisfies either $\partial
U'\subseteq K_{1}$ or $\partial U'\subseteq K_{2}$.

Let $\{U_{i}\}_{i\in I}$ be the set of connected components of $U$
whose boundary is contained in $K_{1}$, and let $\{U_{j}\}_{j\in J}$
be the set of connected components of $U$ whose boundary is contained
in $K_{2}$. Let us set
$$U_{1}:=\bigcup_{i\in I}U_{i},
\hspace{4em}
U_{2}:=\bigcup_{j\in J}U_{j}.$$

It is clear that $U_{1}\cap U_{2}=\emptyset$ and $U_{1}\cup U_{2}=U$.
Moreover from statement~(\ref{lemma:bdr-union}) of
Lemma~\ref{lemma:top} we have that
$$\partial U_{1}=\partial\left(\bigcup_{i\in I}U_{i}\right)\subseteq
\clos\left(\bigcup_{i\in I}\partial
U_{i}\right)\subseteq\clos(K_{1}),$$
and similarly for $U_{2}$.
\qed

A first consequence of Lemma~\ref{lemma:absorb} is the following
result, which is the main tool in the proof of
Theorem~\ref{thm:dvt-maps-homot}.

\begin{lem}\label{lemma:five}
     It is not possible to decompose the unit square $[0,1]\times[0,1]$
     as the disjoint union of subsets $\mathcal{A}_{i}$ ($i=0,1,2,3,4$)
     satisfying the following properties:
     \begin{enumerate}
         \renewcommand{\labelenumi}{(A\arabic{enumi})}
         \item  \label{A-edge} $(0,1)\in\mathcal{A}_{2}$ and
         $(1,1)\in\mathcal{A}_{1}$;

         \item  \label{A-lower}
         $\mathcal{A}_{4}$ does not intersect the side
         $[0,1]\times\{1\}$;

         \item  \label{A-parab} $\mathcal{A}_{0}$ does not intersect
         the other three sides;

         \item \label{A-closed}  $\mathcal{A}_{3}$,
         $\mathcal{A}_{2}\cup\mathcal{A}_{4}$, $\mathcal{A}_{1}\cup
                 \mathcal{A}_{3}\cup\mathcal{A}_{4}$ are closed sets;

                 \item \label{A-open} $\mathcal{A}_{0}$ is an open set and
                 $\partial\mathcal{A}_{0}\subseteq\mathcal{A}_{2}\cup
                 \mathcal{A}_{3}\cup\mathcal{A}_{4}$.
     \end{enumerate}
\end{lem}

\prf
Let us set for simplicity $Q:=[0,1]\times[0,1]$. First of all we show
that, up to modifying the sets
$\mathcal{A}_{0}$, \ldots, $\mathcal{A}_{4}$, we can assume that they
fulfil (A\ref{A-edge}) through (A\ref{A-open}) and also the following
additional property:
\begin{enumerate}
         \item[(A6)] $Q\setminus\mathcal{A}_{0}$ is connected.
\end{enumerate}

Let indeed $P$ be the union of the three sides considered in
(A\ref{A-parab}).  By (A\ref{A-parab}) the closed set
$Q\setminus\mathcal{A}_{0}$ contains the connected set $P$.  Let
$\mathcal{V}$ be the connected component of
$Q\setminus\mathcal{A}_{0}$ containing $P$.
Let us set $\widetilde{\mathcal{A}}_{0}:=Q\setminus\mathcal{V}$ and
$\widetilde{\mathcal{A}}_{i}:=\mathcal{A}_{i}\cap\mathcal{V}$ for
$i=1,2,3,4$.

It is easy to see that the sets $\widetilde{\mathcal{A}}_{0}$, \ldots,
$\widetilde{\mathcal{A}}_{4}$ are disjoint and satisfy assumptions
(A\ref{A-edge}) through (A\ref{A-closed}), and (A6).  Moreover
$\widetilde{\mathcal{A}}_{0}$ is open because $\mathcal{V}$ is closed.
Finally, from statement (\ref{lemma:b-b-c}) of Lemma \ref{lemma:top}
we have that
$$\partial\widetilde{\mathcal{A}}_{0}=\partial\mathcal{V}=
\partial(Q\setminus\mathcal{A}_{0})\cap\mathcal{V}=
\partial\mathcal{A}_{0}\cap\mathcal{V}\subseteq (\mathcal{A}_{2}\cup
\mathcal{A}_{3}\cup\mathcal{A}_{4})\cap\mathcal{V}=
\widetilde{\mathcal{A}}_{2}\cup
\widetilde{\mathcal{A}}_{3}\cup\widetilde{\mathcal{A}}_{4},$$
which proves also (A\ref{A-open}).

Roughly speaking, what we have done in this first part of the proof is
to fill the holes of $\mathcal{A}_{0}$ which do not touch $P$, as
shown in the following picture ($P$ is the union of the lower
and lateral sides of the squares).

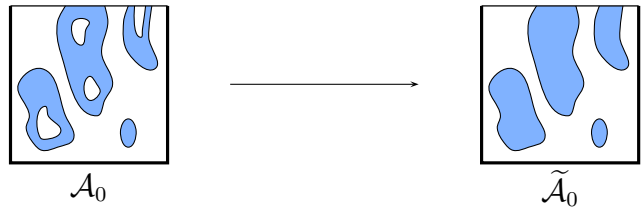
\begin{center}
\psset{unit=1.15ex}
     \pspicture(0,-2)(40,10.5)
     \psset{linewidth=0.5\pslinewidth}
     \newrgbcolor{azero}{0.5 0.7 1}

     \rput{180}(10,10){
 
\pscurve[fillstyle=solid,fillcolor=azero](1,0)(1.2,2)(0.8,4)(2.5,3)(2.8,1)(2.5,0)
     \pscurve[fillstyle=solid,fillcolor=azero]
     (4,0)(3.8,1)(4.2,3)(3.5,5)(4,6)(4.5,7)(6.8,3)(7,1)(6.5,0)
     \psccurve[fillstyle=solid,fillcolor=azero](2,8)(2.5,9)(3,8)(2.5,7.2)
     \psccurve[fillstyle=solid,fillcolor=azero]
     (7,9)(8,9.2)(8.5,8.8)(9,7.8)(8.8,7)(9.5,5)(8,4)(6.5,7)(6,8.5)
     \pscurve[fillstyle=solid,fillcolor=white](1.4,0)(1.6,1.5)(1.7,2)(2.1,0)
     \psccurve[fillstyle=solid,fillcolor=white](5,1)(5.2,1.5)(5.5,2.5)(6,1)
     \psccurve[fillstyle=solid,fillcolor=white](4.5,5)(4.3,5.5)(5,6)(5,4.5)
     \psccurve[fillstyle=solid,fillcolor=white]
     (7.8,8.5)(6.8,8)(7.6,7)(7.8,6.4)(8.1,8.4)
     \psline[linewidth=1.2\pslinewidth](0,0)(10,0)
     \psline[linewidth=3\pslinewidth](0,0)(0,10)(10,10)(10,0)}
     \psline{->}(14,5)(26,5)
     \rput(5,-1.6){$\mathcal{A}_{0}$}

     \psset{origin={20,0}}

     \rput{180}(60,10){
 
\pscurve[fillstyle=solid,fillcolor=azero](1,0)(1.2,2)(0.8,4)(2.5,3)(2.8,1)(2.5,0)
     \pscurve[fillstyle=solid,fillcolor=azero]
     (4,0)(3.8,1)(4.2,3)(3.5,5)(4,6)(4.5,7)(6.8,3)(7,1)(6.5,0)
     \psccurve[fillstyle=solid,fillcolor=azero](2,8)(2.5,9)(3,8)(2.5,7.2)
     \psccurve[fillstyle=solid,fillcolor=azero]
     (7,9)(8,9.2)(8.5,8.8)(9,7.8)(8.8,7)(9.5,5)(8,4)(6.5,7)(6,8.5)
     \psline[linewidth=1.2\pslinewidth](0,0)(10,0)
     \psline[linewidth=3\pslinewidth](0,0)(0,10)(10,10)(10,0)}
     \rput(35,-1.6){$\widetilde{\mathcal{A}}_{0}$}

     \endpspicture
\end{center}

 From now on we drop tildes and we assume that
$\mathcal{A}_{0},\ldots,\mathcal{A}_{4}$ satisfy (A\ref{A-edge})
through (A6).

Since of course $\CH^{1}(Q)=0$, we can apply Lemma
\ref{lemma:absorb} with $Y=Q$, $K_{1}=\mathcal{A}_{3}$,
$K_{2}=\mathcal{A}_{2}\cup \mathcal{A}_{4}$, $U=\mathcal{A}_{0}$.  We
obtain that $\mathcal{A}_{0}$ is the disjoint union of two sets
$\mathcal{A}_{0}'$ and $\mathcal{A}_{0}''$ such that $\partial
\mathcal{A}_{0}'\subseteq \mathcal{A}_{3}$ and $\partial
\mathcal{A}_{0}''\subseteq \mathcal{A}_{2}\cup \mathcal{A}_{4}$.
Together with (A\ref{A-closed}) this implies in particular that
$\mathcal{A}_{0}'\cup \mathcal{A}_{1}\cup \mathcal{A}_{3}\cup
\mathcal{A}_{4}$ and $\mathcal{A}_{0}''\cup \mathcal{A}_{2}\cup
\mathcal{A}_{4}$ are closed subsets of $Q$.

Let us consider now the side $S:=[0,1]\times\{1\}$, which
can be written in the form $$S=\left[S\cap(\mathcal{A}_{0}'\cup
\mathcal{A}_{1}\cup \mathcal{A}_{3}\cup \mathcal{A}_{4})\right]
\cup\left[S\cap(\mathcal{A}_{0}''\cup \mathcal{A}_{2}\cup
\mathcal{A}_{4})\right]=:S_{1}\cup S_{2}.$$

By (A\ref{A-lower}) we have that $S\cap \mathcal{A}_{4}=\emptyset$,
which proves that $S_{1}\cap S_{2}=\emptyset$. By (A\ref{A-edge}) we
have that $(1,1)\in S_{1}$ and $(0,1)\in S_{2}$. Since $S_{1}$ and
$S_{2}$ are closed sets, this contradicts the connectedness of $S$.
\qed

The last lemma is the set-valued extension of a well known result for
continuous functions.

\begin{lem}\label{lemma:corr}
         Let $X$ be a topological space, let $C\subseteq X$ be a closed
         subset, and let $f:C\to\mathcal{P}_{\star}(X)$.  Given $A\subseteq
         C$, let $f(A)$ be the image of $A$, defined as the union of $f(x)$
         when $x$ ranges in $A$.

         If $f$ satisfies (\usc) and (\conn), and $A$ is connected, then
         $f(A)$ is connected.
\end{lem}

\prf
We argue by contradiction.  Let us assume that $U$ and $V$ are open
subsets of $X$ such that $f(A)\cap U$ and $f(A)\cap V$ are nonempty
disjoint sets whose union is $f(A)$.  Let $x\in A$.  Since $f(x)$ is
connected and contained in $f(A)$, it is clear that either
$f(x)\subseteq U$ or $f(x)\subseteq V$.  Therefore if we now define
$$U_{1}:=\{x\in A:\ f(x)\subseteq U\}, \hspace{4em}
V_{1}:=\{x\in A:\ f(x)\subseteq V\},$$
we have found two nonempty disjoint open subsets of $A$ whose union is
$A$.  This contradicts the connectedness of $A$.  \qed

\subsection{Proof of Proposition \ref{prop:dvt-maps}}

As a general fact we recall that, since each $C_{i}$ is a closed set,
the closure $\clos_{i}(Z)$ in $C_{i}$ of any subset $Z\subseteq C_{i}$
coincides with the closure $\clos(Z)$ of $Z$ in $X$.

\paragraph{Statement (1)}

This can be easily proved by induction using the definition of $C_{n}$
and the continuity of $f$.

\paragraph{Statement (\ref{bndr-bndr})}

Since $C_{n+2}$ is closed we have that
$f(\partial_{n+1}C_{n+2})\subseteq f(C_{n+2})\subseteq C_{n+1}$.
Moreover
$$\partial_{n+1}C_{n+2}=\partial_{n+1}(C_{n+1}\setminus
C_{n+2})=\partial_{n+1}A_{n+1}\subseteq\clos_{n+1}(A_{n+1})=\clos(A_{n+1}),$$
hence
$$f(\partial_{n+1}C_{n+2})\subseteq f(\clos(A_{n+1}))\subseteq
\clos(f(A_{n+1}))\subseteq\clos(A_{n})=
\clos_{n}(C_{n}\setminus C_{n+1}).$$

We have thus established that $f(\partial_{n+1}C_{n+2})\subseteq
C_{n+1}\cap\clos_{n}(C_{n}\setminus C_{n+1})$, which is equivalent to
say that $f(\partial_{n+1}C_{n+2})\subseteq\partial_{n}C_{n+1}$.

\paragraph{Statement (\ref{bndr-inside})}

Let us argue by induction.  The case $n=0$ follows from the assumption
that $f(\partial C)\subseteq C$.  Assume now that
$\partial_{n}C_{n+1}\subseteq C_{n+2}$ for some given $n$.  By
statement~(\ref{bndr-bndr}) and the inductive hypothesis we have that
$f(\partial_{n+1}C_{n+2})\subseteq\partial_{n}C_{n+1} \subseteq
C_{n+2}$, which proves that $\partial_{n+1}C_{n+2}\subseteq C_{n+3}$
and completes the induction.

\paragraph{Statement (\ref{ak-closed})}

By statement (\ref{bndr-inside}) we have that
$$\clos(A_{n}) = \clos_{n}(A_{n}) =
A_{n}\cup\partial_{n}A_{n} = A_{n}\cup\partial_{n}(C_{n}\setminus
A_{n}) =
A_{n}\cup\partial_{n}C_{n+1} \subseteq  A_{n}\cup C_{n+2},$$
hence, since $C_{n+2}$ is closed, $\clos(A_{n}\cup
C_{n+2})=\clos(A_{n})\cup\clos(C_{n+2}) \subseteq A_{n}\cup C_{n+2}$,
which completes the proof.
\qed

\subsection{Proof of Theorem \ref{thm:dvt-maps}}

\paragraph{Statement (\ref{stm:dvt2})}
Trivial because $\partial C\subseteq C_{2}$.

\paragraph{Statement (\ref{stm:dvt3})}

If $C=X$, then $\iter(f,C,X)=+\infty$.  If $C$ is a proper subset of
the connected space $X$, then $\partial C\neq\emptyset$, which proves
that $\iter(f,C,X)\geq 2$.  Assume by contradiction that it is exactly
2.  This means that $X=A_{0}\cup A_{1}\cup A_{2}$.  Applying statement
(\ref{ak-closed}) of Proposition~\ref{prop:dvt-maps} with $n=0$ and
$n=1$, we deduce that both $A_{0}\cup A_{2}$ and $A_{1}$ are nonempty
closed sets and this contradicts the connectedness of $X$.

\paragraph{Statement (\ref{stm:dvt4})}

By the previous statement we know that $\iter(f,C,X)\geq 3$.  Assume
by contradiction that it is exactly 3.  This means that $C=A_{1}\cup
A_{2}\cup A_{3}$.  Applying statement (\ref{ak-closed}) of
Proposition~\ref{prop:dvt-maps} with $n=1$ and $n=2$, we deduce that
both $A_{1}\cup A_{3}$ and $A_{2}$ are nonempty closed sets and this
contradicts the connectedness of $C$.

\paragraph{Statement (\ref{stm:dvt6})}

Since $X$ and $C$ are connected, from statement (\ref{stm:dvt4}) we
know that $\iter(f,C,X)\geq 4$.  Assume now by contradiction that it
is exactly 4.  Applying statement~(\ref{ak-closed}) of
Proposition~\ref{prop:dvt-maps} with $n=1,2,3$ we have that $A_{1}\cup
A_{3}\cup A_{4}$, $A_{2}\cup A_{4}$, and $A_{3}$ are nonempty closed
subsets of $X$.

Since $\partial A_{0}=\partial C\subseteq A_{2}\cup A_{3}\cup A_{4}$,
we can apply Lemma \ref{lemma:absorb} with $Y=X$,
$K_{1}=A_{3}$, $K_{2}=A_{2}\cup A_{4}$, $U=A_{0}$. We obtain that
$A_{0}$ is the disjoint union of two sets $A_{0}'$ and $A_{0}''$ such
that $\partial A_{0}'\subseteq A_{3}$ and $\partial A_{0}''\subseteq
A_{2}\cup A_{4}$.  This implies in particular that $A_{0}'\cup
A_{1}\cup A_{3}\cup A_{4}$ and $A_{0}''\cup A_{2}\cup A_{4}$ are
closed subsets of $X$.

Let us consider now the connected set $f(C)$, and let us write
$$f(C)=\left[f(C)\cap(A_{0}'\cup A_{1}\cup A_{3}\cup A_{4})\right]
\cup\left[f(C)\cap(A_{0}''\cup A_{2}\cup A_{4})\right]=:F_{1}\cup F_{2}.$$

Then $F_{1}$ and $F_{2}$ are closed subsets of $f(C)$.  They are also
nonempty because $f(C)$ intersects $A_{1}$, $A_{2}$ and $A_{3}$.
Finally, they are disjoint because $f(C)\cap A_{4}=\emptyset$.  This
contradicts the connectedness of $f(C)$.

\paragraph{Statement (\ref{stm:dvt-infty})}
Let $r:X\setminus\Int{C}{}\to\partial C$ be a retraction, and let
$$g(x):=\left\{
\begin{array}{ll}
     f(x) & \mbox{if }f(x)\in C  \\
     r(f(x)) & \mbox{if }f(x)\not\in\Int{C}{}.
\end{array}
\right.$$

It is not difficult to see that $g:C\to C$ is continuous (one only
needs to verify that it is well defined when $f(x)\in\partial C$).
Since $C$ satisfies the fixed point property there exists $x_{0}\in C$
such that $g(x_{0})=x_{0}$.  We claim that $x_{0}$ is indeed a fixed
point of $f$.

If $f(x_{0})\in C$ then $x_{0}=g(x_{0})=f(x_{0})$ and so $x_{0}$ is
also a fixed point of $f$.  Assume now by contradiction that
$f(x_{0})\not\in C$.  Since $f(\partial C)\subseteq C$, this implies
that $x_{0}\not\in\partial C$.  On the other hand, in this case
$g(x_{0})=r(f(x_{0}))\in\partial C$, which is absurd.  This completes
the proof.  \qed

\subsection{Proof of Theorem \ref{thm:dvt-maps-homot}}

Let $C'$ be a connected component of $C$, and let $X'$ be the
connected component of $X$ containing $C'$.  Since $f$ is homotopic to
the identity it is easy to see that $f$ maps $C'$ to $X'$.  From now
on we can therefore assume that $C$ and $X$ are connected, so that by
statement (\ref{stm:dvt4}) of Theorem \ref{thm:dvt-maps} we have that
$\iter(f,C,X)\geq 4$.

Assume now that it is exactly 4.
Applying statement~(\ref{ak-closed}) of
Proposition~\ref{prop:dvt-maps} with $n=1,2,3$ we have that $A_{1}\cup
A_{3}\cup A_{4}$, $A_{2}\cup A_{4}$, and $A_{3}$ are nonempty closed
subsets of $X$.  Moreover $A_{0}$ is open and $\partial A_{0}\subseteq
A_{2}\cup A_{3}\cup A_{4}$.
\medskip

\emph{Step 1}.  We prove that $\partial C\cap A_{3}\neq\emptyset$.

Let us assume indeed by contradiction that $\partial C=\partial
A_{0}\subseteq A_{2}\cup A_{4}$, hence in particular that $A_{0}\cup
A_{2}\cup A_{4}$ is a closed set.  Now we consider the connected set
$f(C)$ and we write $$f(C)=[f(C)\cap(A_{1}\cup A_{3}\cup A_{4})]\cup
[f(C)\cap(A_{0}\cup A_{2}\cup A_{4})]=:F_{1}\cup F_{2}.$$

Then $F_{1}$ and $F_{2}$ are closed subsets of $f(C)$. They are also
nonempty because $f(C)$ intersects $A_{1}$, $A_{2}$ and $A_{3}$.
Finally, they are disjoint because $f(C)\cap A_{4}=\emptyset$. This
contradicts the connectedness of $f(C)$.
\medskip

\emph{Step 2}.  Let $x_{0}\in \partial C\cap A_{3}$.  We show that
there exists a continuous curve $\gamma:[0,1]\to C$ such that
$\gamma(0)=x_{0}$ and $\gamma(1)\in\partial C\cap A_{2}$.

To begin with, let us consider the
curve $\gamma_{1}:[0,1]\to C$ defined by $\gamma_{1}(t)=\Phi(x_{0},t)$.
This curve takes its values in $C$ because the homotopy sends
$\partial C$ back to $C$.
We can therefore extend it to a curve $\gamma_{2}:[0,2]\to X$ by
setting
$$\gamma_{2}(t):=\left\{
\begin{array}{ll}
     \gamma_{1}(t) & \mbox{if }t\in[0,1],  \\
     \noalign{\vspace{1 ex}}
     f(\gamma_{1}(t-1)) & \mbox{if }t\in[1,2].
\end{array}
\right.$$

The curve $\gamma_{2}$ is continuous (one only needs to check that it
is well defined for $t=1$).  Moreover $\gamma_{2}(1)=f(x_{0})\in
A_{2}$, $\gamma_{2}(2)=f(f(x_{0}))\in A_{1}$, and for every
$t\in[1,2]$ we have that $\gamma_{2}(t)\in f(C)\subseteq A_{0}\cup
A_{1}\cup A_{2}\cup A_{3}$.

We claim that $\gamma_{2}(t)\in A_{0}$ for some $t\in[1,2]$.  Assume
indeed that $\gamma_{2}(t)\in A_{1}\cup A_{2}\cup A_{3}$ for every
$t\in[1,2]$.  Then
$$[1,2] =  \{t\in[1,2]:\ \gamma_{1}(t)\in A_{2}\cup
A_{4}\}\cup\{t\in[1,2]:\ \gamma_{2}(t)\in A_{1}\cup A_{3}\cup
A_{4}\}=:I_{1}\cup I_{2}.$$

Thus $I_{1}$ and $I_{2}$ are closed sets, and they are nonempty
because $1\in I_{1}$ and $2\in I_{2}$. Moreover they are disjoint
because $\gamma_{2}([1,2])\cap A_{4}=\emptyset$, and this contradicts
the connectedness of $[1,2]$.

Let us set now $t_{\star}:=\inf\{t\in[1,2]:\ \gamma_{2}(t)\in
A_{0}\}$.  From the definition of infimum it is clear that
$\gamma_{2}(t)\in C$ for every $t\in[0,t_{\star}]$ and
$\gamma_{2}(t_{\star})\in\partial A_{0}=\partial C$.  We claim that
$\gamma_{2}(t_{\star})\in A_{2}$.  Let us consider indeed
$$[1,t_{\star}]=\{t\in[1,t_{\star}]:\ \gamma_{2}(t)\in A_{1}\cup
A_{3}\}\cup
\{t\in[1,t_{\star}]:\ \gamma_{2}(t)\in A_{2}\}.$$

Once again the two sets in the right hand side are closed and
disjoint, and the second one is nonempty because it contains $t=1$.
By the connectedness of $[1,t_{\star}]$ it follows that the first one
is empty and therefore $\gamma_{2}(t_{\star})\in A_{2}$.

The curve $\gamma$ we are looking for is just (a reparametrization of)
the restriction of $\gamma_{2}$ to the interval $[0,t_{\star}]$.
\medskip

\emph{Step 3}. Let $\gamma$ be the curve of step 2, and let
$$\mathcal{A}_{i}:=\{(\tau,t)\in[0,1]\times[0,1]:\
\Phi(\gamma(\tau),t)\in A_{i}\}$$
for $i=0,1,2,3,4$.  If we show that the $\mathcal{A}_{i}$'s satisfy
assumptions (A\ref{A-edge}) through (A\ref{A-open}) of Lemma
\ref{lemma:five} we have a contradiction.

Since $\gamma(0)\in A_{3}$ we have that
$\Phi(\gamma(0),1)=f(\gamma(0))\in A_{2}$, hence
$(0,1)\in\mathcal{A}_{2}$.  Since $\gamma(1)\in A_{2}$ we have that
$\Phi(\gamma(1),1)=f(\gamma(1))\in A_{1}$, hence
$(1,1)\in\mathcal{A}_{1}$.  This proves (A\ref{A-edge}).

Since the image of $f$ is contained in $A_{0}\cup A_{1}\cup A_{2}\cup
A_{3}$ it follows that $\Phi(\gamma(\tau),1)=f(\gamma(\tau))\not\in
A_{4}$ for every $\tau\in[0,1]$, which proves (A\ref{A-lower}).

Since $\gamma(0)$ and $\gamma(1)$ belong to $\partial C$, and $\Phi$
sends $\partial C$ back to $C$, we have that $\Phi(\gamma(0),t)$ and
$\Phi(\gamma(1),t)$ are in $C$ for every $t\in [0,1]$.  Since also
$\Phi(\gamma(\tau),0)=\gamma(\tau)\in C$ for every $\tau\in[0,1]$,
this proves (A\ref{A-parab}).

Finally, (A\ref{A-closed}) and (A\ref{A-open}) follow from the
continuity of $\Phi(\gamma(\tau),t)$ and the analogous properties of
the $A_{i}$'s.
\qed

\subsection{Proof of Proposition \ref{prop:dvt-corr}}

\paragraph{Statement (\ref{stm:c-cnclosed})}
This can be easily proved by induction using the definition of $C_{n}$
and the upper semicontinuity of $f$.

\paragraph{Statement (\ref{stm:c-a01open})}
The set $A_{0}=X\setminus C$ is open because $C$ is closed. Now since
$$A_{1}=\{x\in C:\ f(x)\cap C=\emptyset\}=
\{x\in C:\ f(x)\subseteq A_{0}\},$$
and since $f$ satisfies (\usc), we have that $A_{1}$ is an open subset of
$C$. In order to conclude that it is also an open subset of $X$ it
suffices to prove that $A_{1}\cap\partial C=\emptyset$. This follows
from (\bw).

\paragraph{Statement (\ref{stm:c-bdr013})}
Let $x\in\partial_{1}C_{2}$.  Since $C_{2}$ is closed we have that
$x\in C_{2}$, hence either $x\in C_{3}$ or $x\in A_{2}$.  Let us
assume by contradiction that $x\in A_{2}$.  Then $f(x)\subseteq
A_{0}\cup A_{1}$ and $f(x)\cap A_{1}\neq\emptyset$.  Since $A_{0}$ and
$A_{1}$ are open sets, and $f(x)$ is connected, we have that
$f(x)\subseteq A_{1}$.  This means that actually $A_{2}=\{x\in C_{1}:\
f(x)\subseteq A_{1}\}$, and thus it is an open subset of $C_{1}$
contained in $C_{2}$.  Therefore if $x\in A_{2}$ then
$x\in\Int{C_{2}}{1}$, which contradicts the initial assumption that
$x\in\partial_{1}C_{2}$.

\paragraph{Statement (\ref{stm:c-a02open})}
The argument is the same used in the proof of statement
(\ref{ak-closed}) of Proposition \ref{prop:dvt-maps}.  Since $C_{1}$
is closed and $\partial_{1}C_{2}\subseteq C_{3}$ we have that
$$\clos(A_{1}) = \mathrm{Clos}_{1}(A_{1}) = A_{1}\cup\partial_{1}A_{1} =
A_{1}\cup\partial_{1}(C_{1}\setminus
A_{1}) = A_{1}\cup\partial_{1}C_{2} \subseteq A_{1}\cup C_{3},$$
hence, since $C_{3}$ is closed, $\clos(A_{1}\cup C_{3})=
\clos(A_{1})\cup\clos(C_{3})=A_{1}\cup C_{3}$.

\paragraph{Statement (\ref{stm:c-bdr124})}
We argue more or less as in the proof of statement
(\ref{stm:c-bdr013}).

Let $x\in\partial_{2}C_{3}$.  Since $C_{3}$ is closed we have that
$x\in C_{3}$, hence either $x\in C_{4}$ or $x\in A_{3}$.  Let us
assume by contradiction that $x\in A_{3}$.  Then $f(x)\subseteq
(A_{0}\cup A_{2})\cup A_{1}$ and $f(x)\cap A_{2}\neq\emptyset$.  Since
$A_{0}\cup A_{2}$ and $A_{1}$ are open sets, and $f(x)$ is connected,
we have that $f(x)\subseteq A_{0}\cup A_{2}$.  This means that
actually $A_{3}=\{x\in C_{2}:\ f(x)\subseteq A_{0}\cup A_{2}\}$, and
thus it is an open subset of $C_{2}$ contained in $C_{3}$.  Therefore
if $x\in A_{3}$ then $x\in\Int{C_{3}}{2}$, which contradicts the
initial assumption that $x\in\partial_{2}C_{3}$.

\paragraph{Statement (\ref{stm:c-a013open})}
Same proof of statement (\ref{stm:c-a02open}) with indices increased by
1.

\paragraph{Statement (\ref{stm:c-fbdr34})}
Let $x\in\partial_{3}C_{4}\cap A_{4}$.  Since $x\in A_{4}$ we know
that $f(x)\cap A_{3}\neq\emptyset$ and $f(x)\subseteq A_{0}\cup
A_{1}\cup A_{2}\cup A_{3}$.  Let us assume by contradiction that
$f(x)\cap A_{2}=\emptyset$, hence that $f(x)$ is contained in the open
set $A_{0}\cup A_{1}\cup A_{3}$.  Now consider $U:=\{x\in C_{3}:\
f(x)\subseteq A_{0}\cup A_{1}\cup A_{3}\}$.  It is an open subset of
$C_{3}$ which is contained in $C_{4}$ (all points in $U$ lie indeed in
$C_{4}$).  Since $x\in U$, we conclude that $x\in\Int{C_{4}}{3}$,
which contradicts the initial assumption that $x\in\partial_{3}C_{4}$.
\qed

\subsection{Proof of Theorem \ref{thm:dvt-corr}}

\paragraph{Statement (\ref{stm:dvtc2})}
Trivial because $\partial C\subseteq C_{2}$.

\paragraph{Statement (\ref{stm:dvtc3})}

If $C=X$, then $\iter(f,C,X)=+\infty$.  If $C$ is a proper subset of
the connected space $X$, then $\partial C\neq\emptyset$, which proves
that $\iter(f,C,X)\geq 2$.  Assume by contradiction that it is exactly
2.  This means that $X=A_{0}\cup A_{1}\cup A_{2}$.  By statements
(\ref{stm:c-a01open}) and (\ref{stm:c-a02open}) of
Proposition~\ref{prop:dvt-corr} we know that both $A_{0}\cup A_{2}$
and $A_{1}$ are nonempty open sets and this contradicts the
connectedness of $X$.

\paragraph{Statement (\ref{stm:dvtc4})}

By the previous statement we know that $\iter(f,C,X)\geq 3$.  Assume
by contradiction that it is exactly 3.  This means that $C=A_{1}\cup
A_{2}\cup A_{3}$.  By statements (\ref{stm:c-a02open}) and
(\ref{stm:c-a013open}) of Proposition~\ref{prop:dvt-corr} we know that
in this case both $A_{1}\cup A_{3}$ and $A_{2}$ are nonempty closed
sets and this contradicts the connectedness of $C$.

\paragraph{Statement (\ref{stm:dvtc5})}

Since $X$ and $C$ are connected, from statement (\ref{stm:dvtc4}) we
know that $\iter(f,C,X)\geq 4$.  Assume now by contradiction that it
is exactly 4.  From statements~(\ref{stm:c-cnclosed}),
(\ref{stm:c-a02open}), and (\ref{stm:c-a013open}) of
Proposition~\ref{prop:dvt-corr} we know that $A_{1}\cup A_{3}\cup
A_{4}$, $A_{2}\cup A_{4}$ and $A_{3}\cup A_{4}=C_{3}$ are closed sets,
but we don't know whether $A_{3}$ is closed or not.

Let us prove that in any case $A_{3}\cap\partial C$ is closed.
Indeed, since
$$\clos(A_{3}) = \mathrm{Clos}_{3}(A_{3}) = A_{3}\cup\partial_{3}A_{3}
= A_{3}\cup\partial_{3}(C_{3}\setminus A_{3}) =
A_{3}\cup\partial_{3}C_{4},$$
and since $\partial_{3}C_{4}\subseteq C_{4}=A_{4}$, we have that
$A_{3}\cap\partial C$ is closed if and only if
$\partial_{3}C_{4}\cap\partial C=\emptyset$.  Let us assume by
contradiction that there exists $x\in\partial_{3}C_{4}\cap\partial C$.
By (\bs) we have that $f(x)\subseteq C$, hence $f(x)\cap
A_{0}=\emptyset$ and therefore
$$f(x)=\left[f(x)\cap(A_{1}\cup
A_{3}\cup A_{4})\right]\cup
\left[f(x)\cap(A_{2}\cup A_{4})\right]=:F_{1}\cup F_{2}.$$

Thus $F_{1}$ and $F_{2}$ are closed subsets of $f(x)$.  Moreover,
since $x\in \partial_{3}C_{4}=\partial_{3}C_{4}\cap A_{4}$, from
statement (\ref{stm:c-fbdr34}) of Proposition \ref{prop:dvt-corr} we
deduce that $F_{1}$ and $F_{2}$ are nonempty.  Finally, they are
disjoint because $f(x)\cap A_{4}=\emptyset$.  This contradicts the
connectedness of $f(x)$.

Once we know that $A_{3}\cap\partial C$ is closed we can proceed as in
the case of functions.
We apply Lemma \ref{lemma:absorb} with $Y=X$, $K_{1}=A_{3}$,
$K_{2}=A_{2}\cup A_{4}$, $U=A_{0}$ and we obtain that $A_{0}$ is the
disjoint union of two sets $A_{0}'$ and $A_{0}''$ such that $\partial
A_{0}'\subseteq \clos(A_{3})\subseteq A_{3}\cup A_{4}$ and $\partial
A_{0}''\subseteq \clos(A_{2}\cup A_{4})=A_{2}\cup A_{4}$.  This
implies in particular that $A_{0}'\cup A_{1}\cup A_{3}\cup A_{4}$ and
$A_{0}''\cup A_{2}\cup A_{4}$ are closed subsets of $X$.

Now we consider $f(C)$, which is a connected set because of
Lemma \ref{lemma:corr}, and we write
$$f(C)=\left[f(C)\cap(A_{0}'\cup A_{1}\cup A_{3}\cup A_{4})\right]
\cup\left[f(C)\cap(A_{0}''\cup A_{2}\cup A_{4})\right].$$

Since $f(C)\cap A_{4}=\emptyset$, the two sets in brackets in the
right hand side are disjoint. They are also nonempty because $f(C)$
intersects $A_{1}$, $A_{2}$ and $A_{3}$. Finally, they are closed
subsets of $f(C)$.

This contradicts the connectedness of $f(C)$.
\qed

\section{Examples}
\label{sec:examples}

The first four examples show that the estimates of $\iter(f,C,X)$
given in the first four statements of Theorem \ref{thm:dvt-maps} are
optimal.

\begin{ex}
     \begin{em}
         Let $X:=\{0\}\cup[2,4]$ with the topology inherited as a
         subset of the real line, let $C:=\{0,4\}$, and let
         $f:C\to X$ be defined by $f(0)=3$ and $f(4)=0$.

                 Then $X$, $C$, and $f$ satisfy the assumptions of Problem
                 \ref{Pbm:map} (in this case indeed $\partial C=\{4\}$), and
                 $\iter(f,C,X)=2$.
    \end{em}
\end{ex}

\begin{ex}
     \begin{em}
                 Let $X:=\R$ with the usual topology, let $C:=\{0\}\cup[2,4]$,
                 and let $f:C\to X$ be defined by $f(x)=(x-2)(x-4)/3$.

                 The function $f$ maps $2$ and $4$ to 0, then it maps $0$
                 inside $(2,4)$, and finally it maps the open interval $(2,4)$
                 outside $C$.

                 Therefore $X$, $C$, and $f$
                 satisfy the assumptions of Problem \ref{Pbm:map}.  Moreover
                 $X$ is connected, $C$ is not connected, $C_{2}=\partial
                 C=\{0,2,4\}$, $C_{3}=\{2,4\}$, and $C_{4}=\emptyset$.  In
                 particular $\iter(f,C,X)=3$.
    \end{em}
\end{ex}

\begin{ex}\label{ex:XC-conn}
     \begin{em}
         Let $X=\{(x,y)\in\R^{2}:\ x^{2}+y^{2}=1\}$ be the circle,
         which we parametrize as usually with the angles in $[0,2\pi]$.
         Let $C:=[2\pi/5,8\pi/5]$ (namely $3/5$ of the circle), and let
         $f:C\to X$ be the counterclockwise rotation by $4\pi/5$
         (namely 2/5 of the way around the circle).

		It turns out that $X$, $C$, and $f$ satisfy the assumptions of
		Problem \ref{Pbm:map} (in this case indeed $\partial C$
		consists of the two points corresponding to $2\pi/5$ and
		$8\pi/5$).  Moreover $X$ and $C$ are connected, and it is not
		difficult to see that $\iter(f,C,X)=4$.  The sets
		$C_{1},\ldots,C_{4}$ are represented in the following
		picture.\medskip

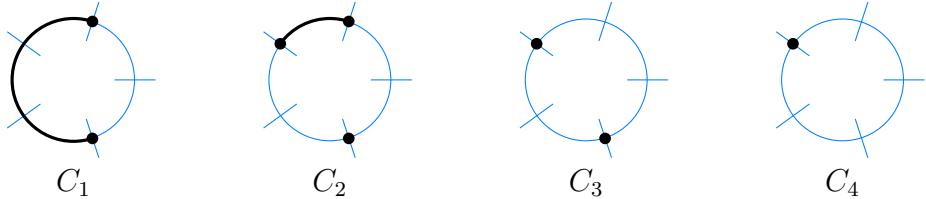
\hfill
\psset{unit=3ex}
\newrgbcolor{Colore}{0 0.5 1}
\pspicture(-2,-3)(2,2)
\SpecialCoor
\psset{linewidth=.5\pslinewidth,linecolor=Colore}
\psline(1;0)(2;0)
\psline(1;72)(2;72)
\psline(1;144)(2;144)
\psline(1;216)(2;216)
\psline(1;288)(2;288)
\pscircle(0,0){1.5}

\psset{linewidth=3\pslinewidth,linecolor=black}
\psarc(0,0){1.5}{72}{288}
\psdots(1.5;72)(1.5;288)
\rput(0,-2.5){$C_{1}$}
\endpspicture
\hfill
\pspicture(-2,-3)(2,2)
\SpecialCoor

\psset{linewidth=.5\pslinewidth,linecolor=Colore}
\psline(1;0)(2;0)
\psline(1;72)(2;72)
\psline(1;144)(2;144)
\psline(1;216)(2;216)
\psline(1;288)(2;288)
\pscircle(0,0){1.5}

\psset{linewidth=3\pslinewidth,linecolor=black}
\psarc(0,0){1.5}{72}{144}
\psdots(1.5;72)(1.5;288)(1.5;144)
\rput(0,-2.5){$C_{2}$}
\endpspicture
\hfill
\pspicture(-2,-3)(2,2)
\SpecialCoor

\psset{linewidth=.5\pslinewidth,linecolor=Colore}
\psline(1;0)(2;0)
\psline(1;72)(2;72)
\psline(1;144)(2;144)
\psline(1;216)(2;216)
\psline(1;288)(2;288)
\pscircle(0,0){1.5}

\psset{linewidth=3\pslinewidth,linecolor=black}
\psdots(1.5;288)(1.5;144)
\rput(0,-2.5){$C_{3}$}
\endpspicture
\hfill
\pspicture(-2,-3)(2,2)
\SpecialCoor

\psset{linewidth=.5\pslinewidth,linecolor=Colore}
\psline(1;0)(2;0)
\psline(1;72)(2;72)
\psline(1;144)(2;144)
\psline(1;216)(2;216)
\psline(1;288)(2;288)
\pscircle(0,0){1.5}
\rput(0,-2.5){$C_{4}$}
\psset{linewidth=3\pslinewidth,linecolor=black}
\psdots(1.5;144)

\endpspicture
\hfill\mbox{}

		We can obviously replace 5 with any greater \emph{odd} integer
		$d$. In this way we obtain a function which can be iterated
		exactly $d-1$ times.

    \end{em}
\end{ex}

\begin{ex}\label{ex:non-cpt}
     \begin{em}
         Let us consider the following subsets of the real plane:
         $$C:=\{(x,0)\in\R^{2}:\ x\in\R\},
         \hspace{3em}
         X_{1}:=\bigcup_{k\in\z}\left([5k,5k+2]\times\R\right).$$

                 Let $X:=X_{1}\cup C$.  Clearly both $X$ and $C$, with the
                 topology inherited as subsets of $\R^{2}$, are connected,
                 simply connected, contractible.  Let $f:C\to X$ be defined by
                 $$f(x,0):=\left\{
         \begin{array}{ll}
             (x+2,0) & \mbox{if $x\in[5k,5k+3]$ for some $k\in\z$,}  \\
             \noalign{\vspace{1 ex}}
             (x+2,|5k+4-x|-1) & \mbox{if $x\in[5k+3,5k+5]$ for some $k\in\z$.}
         \end{array}
         \right.$$

                 Roughly speaking, $C$ is the $x$ axis, $X$ is the union of $C$
                 and some periodically arranged vertical stripes, $f$ is a
                 translation by 2 in the $x$ direction followed by a vertical
                 bending inside the stripes.  The following picture shows the
                 action of $f$ on some points of $C$.

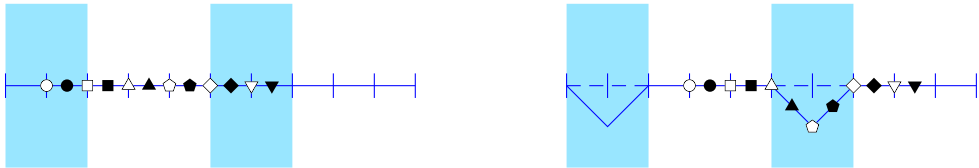
\hfill
\psset{xunit=3ex,yunit=3ex}
\pspicture(0,-3)(10,2.5)
\newrgbcolor{azzurrino}{0.6 0.9 1}
\psframe*[linecolor=azzurrino](0,-2)(2,2)
\psframe*[linecolor=azzurrino](5,-2)(7,2)
\psline[linewidth=.5\pslinewidth,linecolor=blue](0,0)(10,0)
\multido{\n=0+1}{11}{
\psline[linewidth=.5\pslinewidth,linecolor=blue](\n,-0.3)(\n,0.3)}
\psset{dotscale=1.3}
\psdots[dotstyle=o](1,0)
\psdots[dotstyle=square](2,0)
\psdots[dotstyle=triangle](3,0)
\psdots[dotstyle=pentagon](4,0)
\psdots[dotstyle=square,dotangle=45](5,0)
\psdots[dotstyle=triangle,dotangle=180](6,0)
\psdots[dotstyle=*](1.5,0)
\psdots[dotstyle=square*](2.5,0)
\psdots[dotstyle=triangle*](3.5,0)
\psdots[dotstyle=pentagon*](4.5,0)
\psdots[dotstyle=square*,dotangle=45](5.5,0)
\psdots[dotstyle=triangle*,dotangle=180](6.5,0)
\endpspicture
\hfill
\hspace{2em}
\pspicture(0,-3)(10,3)
\newrgbcolor{azzurrino}{0.6 0.9 1}
\psframe*[linecolor=azzurrino](0,-2)(2,2)
\psframe*[linecolor=azzurrino](5,-2)(7,2)
\psline[linewidth=.5\pslinewidth,linecolor=blue](0,0)(1,-1)(2,0)(5,0)(6,-1)(7,0)(10,0)
\psline[linewidth=.5\pslinewidth,linecolor=blue,linestyle=dashed](0,0)(2,0)
\psline[linewidth=.5\pslinewidth,linecolor=blue,linestyle=dashed](5,0)(7,0)
\multido{\n=0+1}{11}{
\psline[linewidth=.5\pslinewidth,linecolor=blue](\n,-0.3)(\n,0.3)}
\psset{dotscale=1.3}
\psdots[dotstyle=o](3,0)
\psdots[dotstyle=square](4,0)
\psdots[dotstyle=triangle](5,0)
\psdots[dotstyle=pentagon](6,-1)
\psdots[dotstyle=square,dotangle=45](7,0)
\psdots[dotstyle=triangle,dotangle=180](8,0)
\psdots[dotstyle=*](3.5,0)
\psdots[dotstyle=square*](4.5,0)
\psdots[dotstyle=triangle*](5.5,-0.5)
\psdots[dotstyle=pentagon*](6.5,-0.5)
\psdots[dotstyle=square*,dotangle=45](7.5,0)
\psdots[dotstyle=triangle*,dotangle=180](8.5,0)
\endpspicture
\hfill
\mbox{}

                 The boundary of $C$ in $X$ is the union of the segments of the
                 form $[5k,5k+2]\times\{0\}$ (the intersection of $C$ with the
                 vertical stripes).  The function $f$ just translates
                 these segments in the $x$ direction, keeping them inside $C$.
                 Therefore all the assumptions of Problem \ref{Pbm:map}
                 are satisfied.  It is not difficult to check that
                 $\iter(f,C,X)=5$, and the sets $C_{2},\ldots,C_{5}$ are those
                 represented in the following picture (we represent only one
                 period, of course).

         \hfill
         \psset{unit=2.5ex}
         
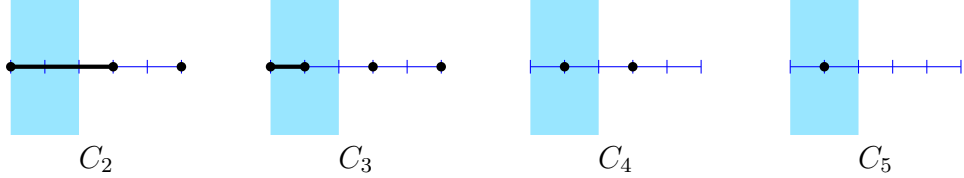
\pspicture(0,-3)(5,3)
         \newrgbcolor{azzurrino}{0.6 0.9 1}
         \psframe*[linecolor=azzurrino](0,-2)(2,2)
         \psline[linewidth=.5\pslinewidth,linecolor=blue](0,0)(5,0)
         \multido{\n=0+1}{6}{
         \psline[linewidth=.5\pslinewidth,linecolor=blue](\n,-0.2)(\n,0.2)}
         \psline[linewidth=2\pslinewidth](0,0)(3,0)
         \psdots(0,0)(3,0)(5,0)
         \rput(2.5,-2.7){$C_{2}$}
         \endpspicture
         \hfill
         \pspicture(0,-3)(5,3)
         \newrgbcolor{azzurrino}{0.6 0.9 1}
         \psframe*[linecolor=azzurrino](0,-2)(2,2)
         \psline[linewidth=.5\pslinewidth,linecolor=blue](0,0)(5,0)
         \multido{\n=0+1}{6}{
         \psline[linewidth=.5\pslinewidth,linecolor=blue](\n,-0.2)(\n,0.2)}
         \psline[linewidth=2\pslinewidth](0,0)(1,0)
         \psdots(0,0)(3,0)(1,0)(5,0)
         \rput(2.5,-2.7){$C_{3}$}
         \endpspicture
         \hfill
         \pspicture(0,-3)(5,3)
         \newrgbcolor{azzurrino}{0.6 0.9 1}
         \psframe*[linecolor=azzurrino](0,-2)(2,2)
         \psline[linewidth=.5\pslinewidth,linecolor=blue](0,0)(5,0)
         \multido{\n=0+1}{6}{
         \psline[linewidth=.5\pslinewidth,linecolor=blue](\n,-0.2)(\n,0.2)}
         \psdots(3,0)(1,0)
         \rput(2.5,-2.7){$C_{4}$}
         \endpspicture
         \hfill
         \pspicture(0,-3)(5,3)
         \newrgbcolor{azzurrino}{0.6 0.9 1}
         \psframe*[linecolor=azzurrino](0,-2)(2,2)
         \psline[linewidth=.5\pslinewidth,linecolor=blue](0,0)(5,0)
         \multido{\n=0+1}{6}{
         \psline[linewidth=.5\pslinewidth,linecolor=blue](\n,-0.2)(\n,0.2)}
         \psdots(1,0)
         \rput(2.5,-2.7){$C_{5}$}
         \endpspicture
         \hfill
         \mbox{}

                 The function $f$ is also homotopic to the identity in the
                 sense of Problem~\ref{Pbm:map-homot}, since both are homotopic
                 to the translation by 2 in the $x$ direction.  Therefore also
                 the assumptions of Problem \ref{Pbm:map-homot} are satisfied,
                 and this shows the optimality of the estimate of
                 $\iter(f,C,X)$ given in Theorem~\ref{thm:dvt-maps-homot}.

     \end{em}
\end{ex}

Note that in Example \ref{ex:non-cpt} above the set $C$ is not
compact.  At the present we have no example of a function $f:C\to X$
satisfying the assumptions of Problem \ref{Pbm:map} with $X$ simply
connected, $C$ compact and connected, and $\iter(f,C,X)=5$.

As we have seen, Example \ref{ex:non-cpt} above shows also the
optimality of Theorem~\ref{thm:dvt-maps-homot}.  We now give another
example, in which the subset $C$ is not only closed, but also compact.

\begin{ex}\label{ex:stargate-bis}
     \begin{em}
         Let us consider polar coordinates $(\rho,\theta)$ in the
         Euclidean plane.  Let
         \begin{eqnarray*}
             X & := & \{ (\rho\cos\theta,\rho\sin\theta)\in\R^{2}:\
             1\leq\rho \leq 3,\ \theta\in[0,2\pi]\}, \\
             Y_1 & := & \{ (\rho\cos\theta,\rho\sin\theta)\in\R^{2}:\
             1\leq\rho \leq 2,\ \theta\in[0,2\pi]\}, \\
             Y_2 & := & \{ (\rho\cos\theta,\rho\sin\theta)\in\R^{2}:\
             2\leq\rho \leq 3,\ 2\pi/5\leq \theta\leq 8\pi/5\}.
         \end{eqnarray*}

         Let
         $C:= Y_1 \cup Y_2$, and let $f:C\to X$ be the function
         represented in polar coordinates by
         $$(\rho,\theta)\to\left(\frac{5}{2},
         \theta+\frac{4\pi}{5}\right).$$

		We claim that $X$, $C$, and $f$ satisfy the assumptions of
		Problem \ref{Pbm:map-homot}.  Indeed, due to our choice of
		$X$, $\partial C$ contains only the arc with $\rho=2$ and
		$\theta\in[-2\pi/5,2\pi/5]$, and the two line segments with
		$\rho\in[2,3]$ and $\theta\in\{-2\pi/5,2\pi/5\}$.  Therefore
		the function $f$ sends $\partial C$ in the points with
		$\rho=5/2$ and $\theta\in[2\pi/5,6\pi/5]$, hence inside $C$.
		As for the required homotopy, roughly speaking it can be
		constructed in three steps: reduction to the level $\rho=2$,
		rotation, reduction to the level $\rho=5/2$.

         After the first iteration all points have radius equal to
         $5/2$, while with regard to the angle we have the identical
         situation of Example \ref{ex:XC-conn}.  It is now simple to
         see that $\iter(f,C,X)=5$ and the sets $C_{1},\ldots,C_{5}$
         are those represented in the following picture.\medskip

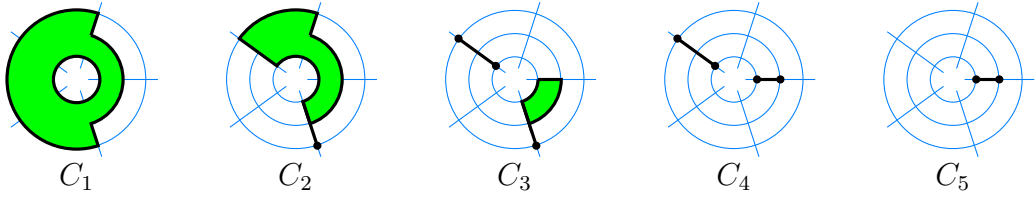
\psset{unit=1.7ex,dotscale=0.7}
\newrgbcolor{Colore}{0 0.5 1}

\hfill
\pspicture(-4,-4.5)(4,4)
\SpecialCoor
\psset{linewidth=.5\pslinewidth,linecolor=Colore}
\psline(0.5;0)(3.5;0)
\psline(0.5;72)(3.5;72)
\psline(0.5;144)(3.5;144)
\psline(0.5;216)(3.5;216)
\psline(0.5;288)(3.5;288)
\pscircle(0,0){1}
\pscircle(0,0){2}
\pscircle(0,0){3}

\psset{linewidth=3\pslinewidth,linecolor=black}
\pscustom[fillstyle=solid,fillcolor=green]{
\psline(2;72)(3;72)
\psarc(0,0){3}{72}{288}
\psline(3;288)(2;288)
\psarc(0,0){2}{288}{360}
\psarcn[liftpen=2](0,0){1}{360}{0}
\psarc[liftpen=2](0,0){2}{0}{72}
}
\rput(0,-4.2)
{$C_{1}$}
\endpspicture
\hfill
\pspicture(-4,-4.5)(4,4)
\SpecialCoor
\psset{linewidth=.5\pslinewidth,linecolor=Colore}
\psline(0.5;0)(3.5;0)
\psline(0.5;72)(3.5;72)
\psline(0.5;144)(3.5;144)
\psline(0.5;216)(3.5;216)
\psline(0.5;288)(3.5;288)
\pscircle(0,0){1}
\pscircle(0,0){2}
\pscircle(0,0){3}

\psset{linewidth=3\pslinewidth,linecolor=black}
\pscustom[fillstyle=solid,fillcolor=green]{
\psline(2;72)(3;72)
\psarc(0,0){3}{72}{144}
\psline(3;144)(1;144)
\psarcn(0,0){1}{144}{288}
\psarc(0,0){2}{288}{72}
}
\psline(1;288)(3;288)
\psdots(3;288)
\rput(0,-4.2)
{$C_{2}$}
\endpspicture
\hfill
\pspicture(-4,-4.5)(4,4)
\SpecialCoor
\psset{linewidth=.5\pslinewidth,linecolor=Colore}
\psline(0.5;0)(3.5;0)
\psline(0.5;72)(3.5;72)
\psline(0.5;144)(3.5;144)
\psline(0.5;216)(3.5;216)
\psline(0.5;288)(3.5;288)
\pscircle(0,0){1}
\pscircle(0,0){2}
\pscircle(0,0){3}

\psset{linewidth=3\pslinewidth,linecolor=black}
\pscustom[fillstyle=solid,fillcolor=green]{
\psline(1;0)(2;0)
\psarcn(0,0){2}{0}{288}
\psline(2;288)(1;288)
\psarc(0,0){1}{288}{0}
}
\psline(1;288)(3;288)
\psline(1;144)(3;144)
\psdots(1;144)(3;144)(3;288)
\rput(0,-4.2)
{$C_{3}$}
\endpspicture
\hfill
\pspicture(-4,-4.5)(4,4)
\SpecialCoor
\psset{linewidth=.5\pslinewidth,linecolor=Colore}
\psline(0.5;0)(3.5;0)
\psline(0.5;72)(3.5;72)
\psline(0.5;144)(3.5;144)
\psline(0.5;216)(3.5;216)
\psline(0.5;288)(3.5;288)
\pscircle(0,0){1}
\pscircle(0,0){2}
\pscircle(0,0){3}

\psset{linewidth=3\pslinewidth,linecolor=black}
\psline(1;0)(2;0)
\psline(1;144)(3;144)
\psdots(1;0)(2;0)
\psdots(1;144)(3;144)
\rput(0,-4.2)
{$C_{4}$}
\endpspicture
\hfill
\pspicture(-4,-4.5)(4,4)
\SpecialCoor
\psset{linewidth=.5\pslinewidth,linecolor=Colore}
\psline(0.5;0)(3.5;0)
\psline(0.5;72)(3.5;72)
\psline(0.5;144)(3.5;144)
\psline(0.5;216)(3.5;216)
\psline(0.5;288)(3.5;288)
\pscircle(0,0){1}
\pscircle(0,0){2}
\pscircle(0,0){3}

\psset{linewidth=3\pslinewidth,linecolor=black}
\psline(1;0)(2;0)
\psdots(1;0)(2;0)
\rput(0,-4.2)
{$C_{5}$}
\endpspicture
\hfill\mbox{}

    \end{em}
\end{ex}

The following Example refers to the Euclidean case of Problem
\ref{Pbm:main}. It is probably the main example of this paper.

\begin{ex}\label{ex:stargate}
     \begin{em}
                 Let $X:=\R^{2}$ be the Euclidean plane, and let $C$ be as in
                 Example \ref{ex:stargate-bis}.  Let $f:C\to X$ be represented
                 in polar coordinates by
                 $$(\rho,\theta)\to\left(\frac{5-|\rho-2|}{2},
         \theta+\frac{4\pi}{5}\right).$$

                 It is clear that $X$ and $C$ are connected, and $X$ is simply
                 connected.  We claim that $X$, $C$, and $f$ satisfy the
                 assumptions of Problem \ref{Pbm:map}.  In this case indeed
                 $\partial C$ contains also the points in $C$ with $\rho=1$ and
                 $\rho=3$, but the image of these points is contained in the
                 level $\rho=2$, hence inside $C$.  Moreover, the function $f$
                 is homotopic to the identity in the sense of Problem
                 \ref{Pbm:map-homot} (as in Example \ref{ex:stargate-bis} the
                 homotopy can be realized through the level $\rho=2$).

                 After two iterations all points have a radius strictly between
                 $2$ and $3$, while with regard to the angle we have the
                 identical situation of Example \ref{ex:XC-conn}.  It is not
                 difficult to see that $\iter(f,C,X)=6$ and the sets
                 $C_{1},\ldots,C_{6}$ are those represented in the following
                 picture.

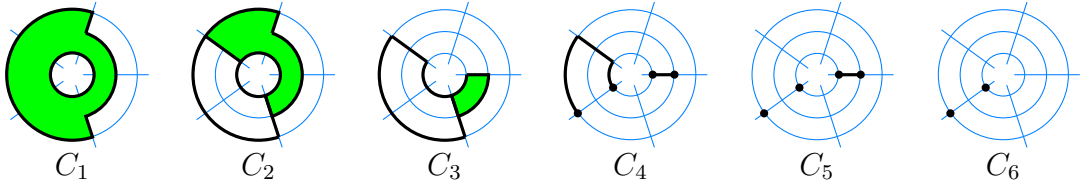
\psset{unit=1.6ex,dotscale=0.7}
\newrgbcolor{Colore}{0 0.5 1}

\hfill
\pspicture(-4,-4.5)(4,4)
\SpecialCoor
\psset{linewidth=.5\pslinewidth,linecolor=Colore}
\psline(0.5;0)(3.5;0)
\psline(0.5;72)(3.5;72)
\psline(0.5;144)(3.5;144)
\psline(0.5;216)(3.5;216)
\psline(0.5;288)(3.5;288)
\pscircle(0,0){1}
\pscircle(0,0){2}
\pscircle(0,0){3}

\psset{linewidth=3\pslinewidth,linecolor=black}
\pscustom[fillstyle=solid,fillcolor=green]{
\psline(2;72)(3;72)
\psarc(0,0){3}{72}{288}
\psline(3;288)(2;288)
\psarc(0,0){2}{288}{360}
\psarcn[liftpen=2](0,0){1}{360}{0}
\psarc[liftpen=2](0,0){2}{0}{72}
}
\rput(0,-4.2)
{$C_{1}$}
\endpspicture
\hfill
\pspicture(-4,-4.5)(4,4)
\SpecialCoor
\psset{linewidth=.5\pslinewidth,linecolor=Colore}
\psline(0.5;0)(3.5;0)
\psline(0.5;72)(3.5;72)
\psline(0.5;144)(3.5;144)
\psline(0.5;216)(3.5;216)
\psline(0.5;288)(3.5;288)
\pscircle(0,0){1}
\pscircle(0,0){2}
\pscircle(0,0){3}

\psset{linewidth=3\pslinewidth,linecolor=black}
\pscustom[fillstyle=solid,fillcolor=green]{
\psline(2;72)(3;72)
\psarc(0,0){3}{72}{144}
\psline(3;144)(1;144)
\psarcn(0,0){1}{144}{288}
\psarc(0,0){2}{288}{72}
}
\psline(1;288)(3;288)
\psarc(0,0){3}{144}{288}
\psarc(0,0){1}{144}{288}
\rput(0,-4.2)
{$C_{2}$}
\endpspicture
\hfill
\pspicture(-4,-4.5)(4,4)
\SpecialCoor
\psset{linewidth=.5\pslinewidth,linecolor=Colore}
\psline(0.5;0)(3.5;0)
\psline(0.5;72)(3.5;72)
\psline(0.5;144)(3.5;144)
\psline(0.5;216)(3.5;216)
\psline(0.5;288)(3.5;288)
\pscircle(0,0){1}
\pscircle(0,0){2}
\pscircle(0,0){3}

\psset{linewidth=3\pslinewidth,linecolor=black}
\pscustom[fillstyle=solid,fillcolor=green]{
\psline(1;0)(2;0)
\psarcn(0,0){2}{0}{288}
\psline(2;288)(1;288)
\psarc(0,0){1}{288}{0}
}
\psline(1;288)(3;288)
\psline(1;144)(3;144)
\psarc(0,0){3}{144}{288}
\psarc(0,0){1}{144}{288}
\rput(0,-4.2)
{$C_{3}$}
\endpspicture
\hfill
\pspicture(-4,-4.5)(4,4)
\SpecialCoor
\psset{linewidth=.5\pslinewidth,linecolor=Colore}
\psline(0.5;0)(3.5;0)
\psline(0.5;72)(3.5;72)
\psline(0.5;144)(3.5;144)
\psline(0.5;216)(3.5;216)
\psline(0.5;288)(3.5;288)
\pscircle(0,0){1}
\pscircle(0,0){2}
\pscircle(0,0){3}

\psset{linewidth=3\pslinewidth,linecolor=black}
\psline(1;0)(2;0)
\psline(1;144)(3;144)
\psarc(0,0){3}{144}{216}
\psarc(0,0){1}{144}{216}
\psdots(1;216)(3;216)
\psdots(1;0)(2;0)
\rput(0,-4.2)
{$C_{4}$}
\endpspicture
\hfill
\pspicture(-4,-4.5)(4,4)
\SpecialCoor
\psset{linewidth=.5\pslinewidth,linecolor=Colore}
\psline(0.5;0)(3.5;0)
\psline(0.5;72)(3.5;72)
\psline(0.5;144)(3.5;144)
\psline(0.5;216)(3.5;216)
\psline(0.5;288)(3.5;288)
\pscircle(0,0){1}
\pscircle(0,0){2}
\pscircle(0,0){3}

\psset{linewidth=3\pslinewidth,linecolor=black}
\psline(1;0)(2;0)
\psdots(1;216)(3;216)
\psdots(1;0)(2;0)
\rput(0,-4.2)
{$C_{5}$}
\endpspicture
\hfill
\pspicture(-4,-4.5)(4,4)
\SpecialCoor
\psset{linewidth=.5\pslinewidth,linecolor=Colore}
\psline(0.5;0)(3.5;0)
\psline(0.5;72)(3.5;72)
\psline(0.5;144)(3.5;144)
\psline(0.5;216)(3.5;216)
\psline(0.5;288)(3.5;288)
\pscircle(0,0){1}
\pscircle(0,0){2}
\pscircle(0,0){3}

\psset{linewidth=3\pslinewidth,linecolor=black}
\psdots(1;216)(3;216)
\rput(0,-4.2)
{$C_{6}$}
\endpspicture
\hfill\mbox{}

\end{em}
\end{ex}

The constructions presented in Example~\ref{ex:XC-conn} and
Example~\ref{ex:stargate} can be extended to higher dimensions,
proving that also the triviality of some higher dimensional
cohomology groups of $C$ and $X$ doesn't imply the existence of an
infinite orbit.

\begin{ex}\label{ex:ndim}
	\begin{em}
		Let $X$ be the unit sphere in $\R^{2n}=\C^{n}$, and let
		$d_{1},\ldots,d_{n}$ be pairwise coprime odd integers each
		greater than or equal to 5. Let us define
		$\lambda_{k}:=\exp(4\pi i/d_{k})$, and
		$$f(z_{1},\ldots,z_{n}):=
		(\lambda_{1}z_{1},\ldots,\lambda_{n}z_{n}).$$

		Let $O$ be the set of points $(z_{1},\ldots,z_{n})\in X$ with
		$0<\arg(z_{k})<4\pi/d_{k}$ for every $k=1,\ldots,n$, and let
		$C:=X\setminus O$.

		For this choice it is clear that $C$ is contractible,
		$H^{i}(X)=0$ for all $1\leq i\leq 2n-1$, and
		$\iter(f,C,X)=d_{1}d_{2}\cdots d_{n}$.

		Again, one can improve the construction via the same process
		as in the passage from Example~\ref{ex:XC-conn} to
		Example~\ref{ex:stargate} to get a similar example with
		$X=\R^{2n}$ and $H^{i}(C)=0$ for all $1\leq i\leq 2n-1$.
	\end{em}
\end{ex}

Now we present two examples concerning set-valued maps. The first one
shows that without connectedness assumptions on the images only
trivial iterations are guaranteed.

\begin{ex}\label{ex:trivial-corr}
     \begin{em}
                 Let $X:=\R$ with the usual topology, and let $C:=[0,4]$.  Let
                 $f:C\to\mathcal{P}_{\star}(X)$ be defined by
                 $$f(x):=\left\{
                 \begin{array}{ll}
                         \{2\} & \mbox{if }x\in[0,1)\cup(3,4],  \\
                         \{2,5\} & \mbox{if }x\in\{1,3\},   \\
                         \{5\} & \mbox{if }x\in(1,3).
                 \end{array}
                 \right.$$

                 It is easy to show that $X$ and $C$ are connected (and even
                 contractible), $f$ satisfies (\usc) and (\bs), and
                 $\iter(f,C,X)=2$.
    \end{em}
\end{ex}

The following Example shows the optimality of statement
(\ref{stm:dvtc4}) in Theorem \ref{thm:dvt-corr}.

\begin{ex}\label{ex:corr-main}
     \begin{em}
         Let $X:=\R^{2}$ with the usual topology, let $Q$ be the
                 square $[0,2]\times[0,2]$, let $S$ be the segment with
                 endpoints $(2,0)$ and $(3,0)$, and let $C:=Q\cup S$. Let
                 $d(x,y)$ denote the distance of the point $(x,y)$ from the
                 boundary of $Q$. Let $f:C\to\mathcal{P}_{\star}(X)$ be defined
                 in the following way:
                 \begin{itemize}
                         \item  if $(x,y)=(2,0)$, then $f(x,y)$ is the 
segment with
                         endpoints $(1,1)$ and $(3,0)$;

                         \item if $(x,y)\in S\setminus Q$, then $f(x,y)$ is the
                         singleton $\{(1,1)\}$;

                         \item if $(x,y)\in Q\setminus S$, then $f(x,y)$ is the
                         singleton $\{(3,d(x,y))\}$.
                 \end{itemize}

                 In a few words, $f$ is single-valued except at $(2,0)$: it
                 send $\partial Q\setminus S$ to $(3,0)$, in turn $(3,0)$ and
                 the rest of $S\setminus Q$ are sent inside the square at
                 $(1,1)$, and the interior of $Q$ is sent outside $C$. Finally,
                 the image of $(2,0)$ is the minimal convex set for which the
                 resulting function turns out to be upper semicontinuous.

                 Therefore $f$ satisfies (\usc), (\bw), (\conn).
                 Moreover $\iter(f,C,X)=4$, and the sets $C_{1}$, \ldots,
                 $C_{4}$ are those represented in the following picture.

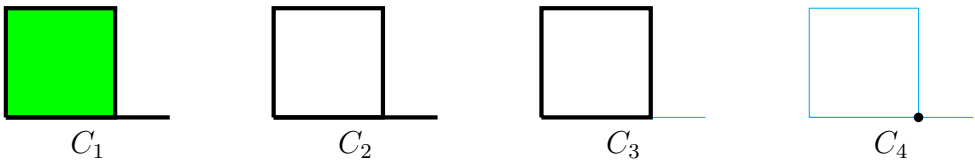
\begin{center}
         \psset{unit=4ex}

         \hfill
         \pspicture(0,-1)(3,2)
 
\psline[linewidth=2\pslinewidth,fillstyle=solid,fillcolor=green](0,0)(2,0)(2,2)(0,2)(0,0)
         \psline[linewidth=2\pslinewidth](0,0)(3,0)
         \rput(1.5,-0.5){$C_{1}$}
         \endpspicture
         \hfill
         \pspicture(0,-1)(3,2)
         \psline[linewidth=2\pslinewidth](0,0)(2,0)(2,2)(0,2)(0,0)
         \psline[linewidth=2\pslinewidth](0,0)(3,0)
         \rput(1.5,-0.5){$C_{2}$}
         \endpspicture
         \hfill
         \pspicture(0,-1)(3,2)
         \psline[linecolor=cyan,linewidth=0.5\pslinewidth](0,0)(3,0)
         \psline[linewidth=2\pslinewidth](0,0)(2,0)(2,2)(0,2)(0,0)
         \rput(1.5,-0.5){$C_{3}$}
         \endpspicture
         \hfill
         \pspicture(0,-1)(3,2)
         \psline[linecolor=cyan,linewidth=0.5\pslinewidth](0,0)(3,0)
 
\psline[linecolor=cyan,linewidth=0.5\pslinewidth](0,0)(2,0)(2,2)(0,2)(0,0)
         \psdot(2,0)
         \rput(1.5,-0.5){$C_{4}$}
         \endpspicture
         \hfill\mbox{}
\end{center}

    \end{em}
\end{ex}

\section{Stochastic games and iterations}\label{sec:games}

\paragraph{Stochastic games}

In literature there are different and sometimes contradictory 
definitions of stochastic game, with different levels of generality. 
For the sake of simplicity, we present here a definition which is 
quite restrictive, but yet enough to describe this field and its open 
problems. In this paragraph we stick to the notations 
of~\cite{v1}. 

A stochastic game is played by a finite set $P$ of players over a
finite set $S$ of possible states, and involves a finite or countable
number of stages.  At each stage $n$ the game is in some state
$s_{n}\in S$.  If this is not the last stage of play, each player
chooses an action in a finite set $A$ of possible options, and the
state of the game changes to some $s_{n+1}$ which is a (possibly
random) function of $s_{n}$ and of the element of $A^{P}$
representing the players' choices.  Moreover, each player receives a
payoff, which also depends on $s_{n}$ and on the actions selected.
Then the game moves to next stage.  In all stages, all players have
complete knowledge of the past history of play, of the present state,
and of the present options and their consequences.  The only
uncertainty concerns what the other players will do in the present and
in the future.

Let $H_{n}:=S\times (S\times A^{P})^{n-1}$ denote the set of possible
histories up to stage $n$, and let $H$ denote the union of $H_{n}$ as
$n$ varies over all stages, namely the set of all finite histories.  A
\emph{strategy} for a given player is a function
$H\to\mathrm{Prob}(A)$, where $\mathrm{Prob}(A)$ is the space of
probability measures on $A$ (actually it is a simplex).  When the game
is in stage $n$, and $h_{n}\in H_{n}$ is the history up to that stage,
then the value of the strategy in $h_{n}$ is the lottery used by the
player in order to select next action.

Up to introducing a cumbersome notation, one could also admit that
states are a set $S_{n}$ depending on the stage $n$, and possible
actions are a set $A_{n,s,p}$ depending on the stage $n$, on the state
$s$, and on the player $p$.  Of course in this more general setting
the definitions of histories and strategies need to be changed
accordingly.

For any $\ep\geq 0$, an {\em $\ep$-equilibrium} in a game is a profile
of strategies, one for each player, such that no player can gain in
expected payoff by more than $\ep$ by choosing a different strategy,
given that all the other players do not change their strategies (for
more precision, see~\cite{v1}).  An \emph{equilibrium} is a
$0$-equilibrium.  We say that approximate equilibria exist if there
exists an $\ep$-equilibrium for every $\ep>0$.

It is well known that equilibria exist whenever the stochastic game
has finitely many stages.  This is a celebrated result by \textsc{J.\
Nash}~\cite{nash}.

When there are infinitely many stages of play, things are more
complex.  In the special case where the number of players is two,
\textsc{N.\ Vieille}~\cite{v1,v2,v3} proved existence of approximate
equilibria.  With three or more players, it is not known whether all
stochastic games have approximate equilibria.  \textsc{R.\ Aumann},
during his Address to the first world congress of ``The Game Theory
Society'' (GAMES 2000, Bilbao), stated that this question is the most
important open problem of mathematical game theory today.

\paragraph{Quitting games}

Quitting games are a special class of stochastic games with a very
simple structure.  At any stage of a quitting game, each player has
only two actions, \texttt{c} for continue and \texttt{q} for quit.  As
soon as one or more of the players at any stage chooses \texttt{q},
the game stops, and players receive payoffs, which depend on the
subset of players that choose simultaneously the action \texttt{q}.
Whenever all players choose \texttt{c}, the game goes to the next
stage, and all players receive the payoff of $0$ for that stage.
Quitting games were studied first by \textsc{J.\ Flesch}, \textsc{F.\
Thuijsman} and \textsc{O.\ J.\ Vrieze} in~\cite{FTV}, but modelled
first in full generality by \textsc{E.\ Solan} and \textsc{N.\
Vieille}~\cite{s-v}.

The complexity of quitting games lies in the potentially large number
of players.  In the case of two players, one can prove existence of
\emph{stationary} $\ep$-equilibria, namely $\ep$-equilibria where
strategies depend on the past history only through the current state.
In~\cite{FTV} a three-player example was shown, where
$\ep$-equilibrium strategies have a nonconstant cyclic structure.
This motivated the study of the three-player case, solved by
\textsc{E. Solan}~\cite{solan} by proving that approximate equilibria
do exist.  With four or more players, the problem is still open.

In a few words, despite the simpler structure, quitting games are
important for both the positive and negative sides of the question of
whether approximate equilibria exist for general stochastic
games.

\paragraph{Quitting games  and iterations of set-valued maps}

An approach connecting quitting games and topological dynamics has
been introduced by \textsc{E.\ Solan} and \textsc{N.\ Vieille}
in~\cite{s-v}.  We sketch the main steps of this approach
following~\cite{s-v}, to which we refer for further details.

The main idea is to break up the game into infinitely many one-shot
games, namely games played in one stage only.

To this end, let $N:=|P|$ be the number of players.  For every fixed
vector $w\in\R^{N}$, let $\Gamma_{w}$ be the one-shot game where the
payoff vector is the same as in the original quitting game if at least
one player chooses \texttt{q}, and it is $w$ otherwise.  The strategy of a
player in $\Gamma_{w}$ is just the probability to choose \texttt{c},
so that $[0,1]^{N}$ is the set of strategy profiles for $\Gamma_{w}$.
Let $f(w,p)\in\R^{N}$ be the expected payoff vector in the game
$\Gamma_{w}$ when all players perform according to some strategy
profile $p=(p_{1},\ldots,p_{N})\in[0,1]^{N}$, and let $q(p):=1-
p_{1}\cdot\ldots\cdot p_{N}$ be the probability that at least one 
player chooses \texttt{q}.

Now let $\rho$ be a large enough constant, depending only upon payoff
vectors, and let $\ep>0$ be small enough.  Let
$E_{\rho\ep}(w)\subseteq[0,1]^{N}$ be the set of $\rho\ep$-equilibria
for the game $\Gamma_{w}$, which is nonempty because of Nash's
theorem, and let us finally define
\begin{equation}
	F_{\ep}(w):=\left\{f(w,p):p\in E_{\rho\ep}(w),\ q(p)\geq \ep\right\}
	\subseteq\R^{N}.
	\label{defn:Fep-sv}
\end{equation}

We have thus a map $F_{\ep}$ from $\R^{N}$ to subsets of
$\R^{N}$, which can be subjected to iteration.  It is not difficult to
see that this set-valued map is upper semicontinuous.  The problem is
that $F_{\ep}(w)$ might be empty for a set of $w$ which is 
eventually reached by any iteration process, and this would 
prevent infinite orbits from existing.

The main achievements of~\cite{s-v} are providing conditions on the
game that ensure that an infinite orbit exists (Proposition 2.2), and
proving that any orbit defines an equilibrium (Proposition 2.4). This 
relates equilibria of the game with dynamics.

\paragraph{Relations with DVT}

What is actually proved in Proposition~2.2 of~\cite{s-v} is that, 
under suitable conditions on the game, there 
exists a compact set $C\subseteq\R^{N}$ such that 
$F_{\ep}(w)\cap C\neq\emptyset$ for every $w\in C$. This guarantees 
for free the existence of an infinite orbit, and produces a class of 
quitting games with approximate equilibria. The existence of such a 
set is however just a sufficient condition.

Further relations between games and dynamics have been investigated
in~\cite{topological,simon}.  In these papers $F_{\ep}(w)$ is defined
as in (\ref{defn:Fep-sv}), but without condition $q(p)\geq\ep$, so
that now $F_{\ep}(w)$ is trivially nonempty for every $w$.  With such
a definition, Theorem~3 in~\cite{simon} states that (but for trivial
cases) approximate equilibria exist if and only if for every small
enough $\ep >0$ there exist (infinite) orbits of this new $F_{\ep}$
with \emph{unbounded total variation} (the total variation of an orbit
is the sum of distances between consecutive terms).  This is a
necessary and sufficient condition.  The requirement on the total
variation, which in the game context implies eventual quitting with
certainty, rules out fixed points, and makes the problem highly
nontrivial.

Theorem~1 of~\cite{topological} states that there exists a compact
connected set $C\subseteq\R^{N}$ such that $F_{0}$, as a map form $C$
to $\mathcal{P}_{\star}(\R^{N})$, has the following topological
properties.  The set of points $w\in C$ such that $w\in F_{0}(w)$
coincides with $\partial C$, and $F_{0}$ is homotopic to the identity
map on $C$ through a homotopy whose intermediate maps keep again all
points of $\partial C$ fixed.  This homotopy condition motivated our
interest in our Problem~\ref{Pbm:map-homot}.

Several tricks can be devised in order to avoid useless stationary or
converging orbits.  One possibility is taking the set $C$ considered
before, and introducing a new set-valued map $G_{\ep}$ whose graph is
obtained from the graph of $F_{\ep}$ by removing an open set
containing the fixed points $(w,w)$ with $w\in\partial C$.  One can
prove that, if $\ep>0$ and the removed open set is small enough, this
construction gives a well defined set-valued map without fixed points,
and such that for every $w\in\partial C$ the image $G_{\ep}(w)$ still
contains some motion back to the set $C$.  Due to the lack of fixed
points, any orbit of $G_{\ep}$ is a non-converging orbit of $F_{\ep}$,
hence an orbit with unbounded total variation.  We have thus reduced
the problem to a situation similar to our Problem~\ref{Pbm:map} or
Problem~\ref{Pbm:corr}.

\noindent
\psset{unit=2.5ex}
\hfill
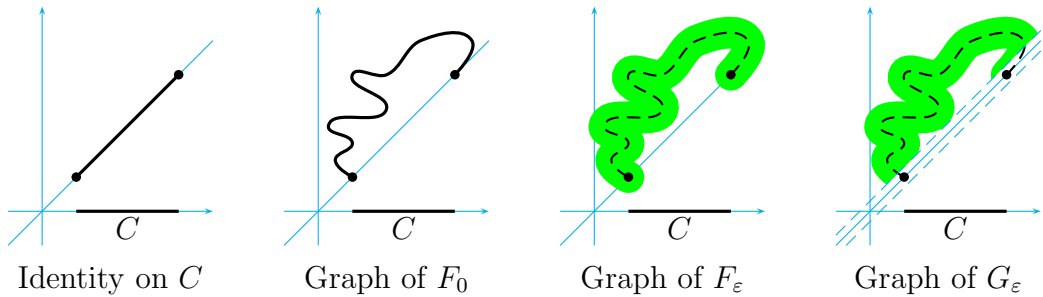
\pspicture(-1,-3)(5,6.5)
\psline[linewidth=0.5\pslinewidth,linecolor=cyan]{->}(-1,0)(5,0)
\psline[linewidth=0.5\pslinewidth,linecolor=cyan]{->}(0,-1)(0,6)
\psline[linewidth=0.5\pslinewidth,linecolor=cyan](-1,-1)(5,5)
\psline[linewidth=1.5\pslinewidth](1,0)(4,0)
\psline[linewidth=1.5\pslinewidth](1,1)(4,4)
\psdots(1,1)(4,4)
\rput(2.5,-0.5){$C$}
\rput(2,-2){Identity on $C$}
\endpspicture
\hfill
\pspicture(-1,-3)(5,6.5)
\psline[linewidth=0.5\pslinewidth,linecolor=cyan]{->}(-1,0)(5,0)
\psline[linewidth=0.5\pslinewidth,linecolor=cyan]{->}(0,-1)(0,6)
\psline[linewidth=0.5\pslinewidth,linecolor=cyan](-1,-1)(5,5)
\psline[linewidth=1.5\pslinewidth](1,0)(4,0)
\pscurve[linewidth=1.5\pslinewidth](1,1)(0.5,1.5)(1,2)(0.3,2.5)(2,3)(1,4)
(2,4)(3,5)(4.5,5)(4,4)
\psdots(1,1)(4,4)
\rput(2.5,-0.5){$C$}
\rput(2,-2){Graph of $F_{0}$}
\endpspicture
\hfill
\pspicture(-1,-3)(5,6.5)
\psline[linewidth=0.5\pslinewidth,linecolor=cyan]{->}(-1,0)(5,0)
\psline[linewidth=0.5\pslinewidth,linecolor=cyan]{->}(0,-1)(0,6)
\psline[linewidth=0.5\pslinewidth,linecolor=cyan](-1,-1)(5,5)
\psline[linewidth=1.5\pslinewidth](1,0)(4,0)
\pscurve[linewidth=15\pslinewidth,linecolor=green]{c-c}(1,1)(0.5,1.5)(1,2)(0.3,2.5)(2,3)(1,4)
(2,4)(3,5)(4.5,5)(4,4)
\pscurve[linewidth=1\pslinewidth,linestyle=dashed](1,1)(0.5,1.5)(1,2)(0.3,2.5)(2,3)(1,4)
(2,4)(3,5)(4.5,5)(4,4)
\psdots(1,1)(4,4)
\rput(2.5,-0.5){$C$}
\rput(2,-2){Graph of $F_{\ep}$}
\endpspicture
\hfill
\pspicture(-1,-3)(5,6.5)
\psline[linewidth=0.5\pslinewidth,linecolor=cyan]{->}(-1,0)(5,0)
\psline[linewidth=0.5\pslinewidth,linecolor=cyan]{->}(0,-1)(0,6)
\psline[linewidth=0.5\pslinewidth,linecolor=cyan](-1,-1)(5,5)
\psline[linewidth=1.5\pslinewidth](1,0)(4,0)
\psclip{\psline[linestyle=none](-1,-0.7)(-1,7)(6.7,7)(-1,-0.7)}
\pscurve[linewidth=15\pslinewidth,linecolor=green]{c-c}(1,1)(0.5,1.5)(1,2)(0.3,2.5)(2,3)(1,4)
(2,4)(3,5)(4.5,5)(4,4)
\endpsclip
\pscurve[linewidth=1\pslinewidth,linestyle=dashed](1,1)(0.5,1.5)(1,2)(0.3,2.5)(2,3)(1,4)
(2,4)(3,5)(4.5,5)(4,4)
\psline[linewidth=0.5\pslinewidth,linecolor=cyan,linestyle=dashed](-0.7,-1)(5,4.7)
\psline[linewidth=0.5\pslinewidth,linecolor=cyan,linestyle=dashed](-1,-0.7)(5,5.3)
\psdots(1,1)(4,4)
\rput(2.5,-0.5){$C$}
\rput(2,-2){Graph of $G_{\ep}$}
\endpspicture
\hfill\mbox{}

The picture above is an attempt to represent this situation, with all
limitations of a two dimensional setting.

Approaches of this kind have failed so far to establish the existence of
approximate equilibria for quitting games, mostly due to the lack of
any corresponding theorems of topological dynamics that demonstrates
the existence of infinite orbits from some topological properties
extracted from quitting games.

Indeed, on the game side we suspect that there are quitting games that
do not have approximate equilibria, as well as on the dynamics side we
suspect that there are many more plausible ``results'' concerning the
existence of infinite orbits that fail to be true.

\section{Open problems}\label{sec:open}
As mentioned in the introduction, the following is probably the main
question in Discrete Viability Theory.
\begin{open}
         Find nontrivial sufficient conditions on $f$, $C$, $X$ in
         order to have that $\iter(f,C,X)=+\infty$.
\end{open}

Here ``nontrivial'' means that these conditions should be satisfied by
reasonable classes of functions $f$ without fixed points and with
$f(C)\not\subseteq C$.

A first step in this direction could be to understand whether
strengthening the topological assumptions on $f$, $C$, $X$ guarantees
further iterations. This leads to the second question.
\begin{open}
         Under the assumptions of Problem~\ref{Pbm:map} find intermediate
         results between statement~(\ref{stm:dvt6}) and
         statement~(\ref{stm:dvt-infty}) of Theorem~\ref{thm:dvt-maps}.
\end{open}

Example \ref{ex:non-cpt} shows that 5 iterations is the most one can
expect even when $C$ and $X$ are contractible, hence as simple as
possible from the topological point of view.  This seems to be the
tombstone on the search of further iterations.  Nevertheless, we point
out once again that in that example $C$ is not compact.  So a new
frontier is understanding the role played by compactness in this
subject, even in the simpler case.
\begin{open}
         Find the maximal number of iterations which are assured under
         the assumptions of Problem~\ref{Pbm:main}.
\end{open}

We know that this number is at least 5 and at most 6.  We also know
that if this number is 6 it is a matter of compactness. If this is the
case, then we can ask ourselves what happens with further topological
requirements on $C$, for example if $\CH^{1}(C)=0$ (just to rule out
Example~\ref{ex:stargate}). We know from Example~\ref{ex:ndim} that we
cannot expect an infinite orbit, but maybe the lower bound on the
number of iterations increases.

Finally, a technical point for topologists.
\begin{open}
         Find the minimal assumptions on $Y$ under which
         Lemma~\ref{lemma:absorb} can be proved.
\end{open}

\paragraph{Acknowledgments}
This research was supported financially by the German Science
Foundation (Deutsche Forschungsgemeinschaft).

\label{NumeroPagine}

\end{document}